\documentclass[11pt]{scrartcl}
\pdfoutput=1
\usepackage{amsmath,mathtools}
\usepackage{url,color}
\usepackage[pdftex,colorlinks]{hyperref}
\definecolor{darkblue}{cmyk}{1,0,0,0.8}
\definecolor{darkred}{cmyk}{0,1,0,0.7}
\hypersetup{anchorcolor=black,
  citecolor=darkblue, filecolor=darkblue,
  menucolor=darkblue,pagecolor=darkblue,urlcolor=darkblue,linkcolor=darkblue}
\usepackage[scaled=0.8]{helvet}
\usepackage{amssymb,amsfonts,mathpazo}
\usepackage{eulervm}
\usepackage[labelfont={small,bf},textfont={small},format=plain]{caption}
\usepackage{color}
\newcommand{\rev}[1]{{#1}}
\usepackage{microtype} 
\usepackage{mythmbox}
\typearea{12}

\newtheorem[S]{corollary}{Corollary}[section]
\newtheorem[S]{theorem}[corollary]{Theorem}
\newtheorem[S]{definition}[corollary]{Definition}
\newtheorem[S]{assumption}[corollary]{Assumption}
\newtheorem[S]{lemma}[corollary]{Lemma}
\newtheorem[S]{proposition}[corollary]{Proposition}
\newtheorem[S]{introthm}{Main Result}

\newcommand{\R}{\mathbb{R}}

\newcommand{\Z}{\mathbb{Z}}

\newcommand{\C}{\mathbb{C}}

\newcommand{\id}{I}

\newcommand{\spec}{\operatorname{spec}}

\renewcommand{\mod}{\operatorname{mod}}

\renewcommand{\i}{{\mathrm{i}}}
\renewcommand{\d}{\mathop{}\!\mathrm{d}}
\renewcommand{\Re}{\operatorname{Re}}

\renewcommand{\phi}{\varphi}
\renewcommand{\epsilon}{\varepsilon}

\newcommand{\eop}[1]{\hfill (\emph{This ends the proof of #1.})\quad$\square$}

\allowdisplaybreaks
\title{Generic stabilizability for time-delayed feedback
  control}

\author{Jan
  Sieber\\  
  College of Engineering, Mathematics and Physical
  Sciences,\\ University of Exeter, EX4 4QF, UK}
\begin{document}
\maketitle\begin{abstract}\noindent
  Time delayed feedback control is one of the most successful methods
  to discover dynamically unstable features of a dynamical system in
  an experiment. This approach feeds back only terms that depend on
  the difference between the current output and the output from a
  fixed time $T$ ago. Thus, any periodic orbit of period $T$ in the
  feedback controlled system is also a periodic orbit of the
  uncontrolled system, independent of any modelling assumptions. 

  It has been an open problem whether this approach can be successful
  in general, that is, under genericity conditions similar to those in
  linear  control theory (controllability), or if there are
  fundamental restrictions to time-delayed feedback control. We show
  that there are no restrictions in principle. This paper proves the
  following: for every periodic orbit satisfying a genericity
  condition slightly stronger than classical linear controllability,
  one can find control gains that stabilise this orbit with extended
  time-delayed feedback control.

  While the paper's techniques are based on linear stability analysis,
  they exploit the specific properties of linearisations near
  autonomous periodic orbits in nonlinear systems, and are, thus,
  mostly relevant for the analysis of nonlinear experiments.
\end{abstract}

\section{Introduction}
\label{sec:intro}
\begin{figure}[ht]
  \centering
    \includegraphics[scale=0.5]{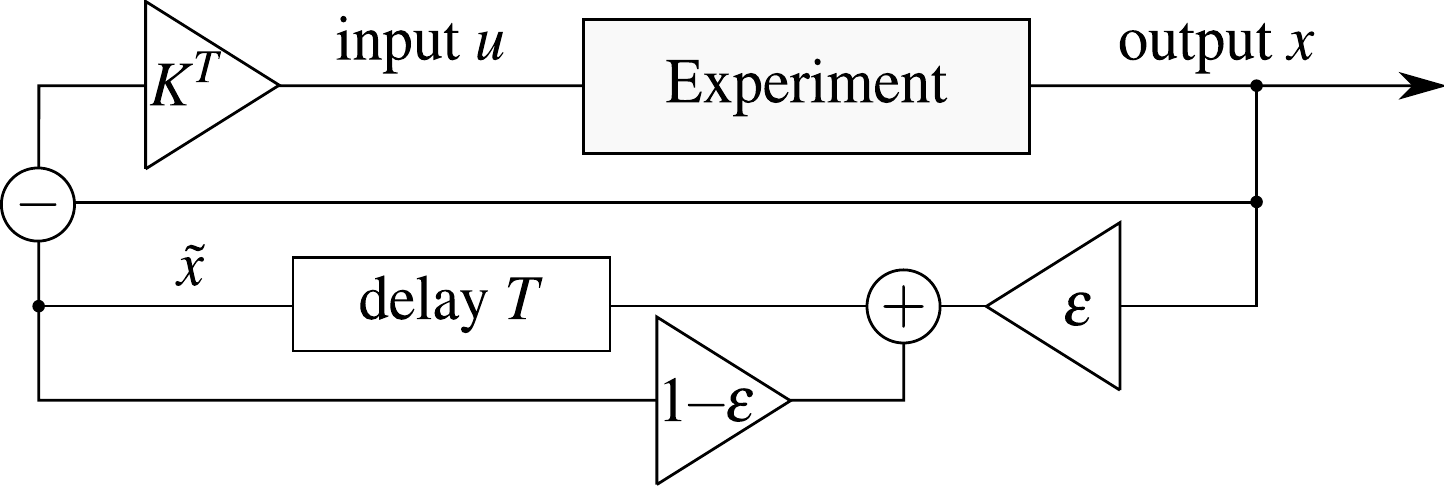}
    \caption{\rev{Block diagram for extended time-delayed feedback
        \eqref{eq:intro:etdf}, as applied to an experiment, for
        example, in \cite{SSG94}.  We prove generic stabilisablity for the
        case when the output $x$ is the whole internal state of the
        dynamical system and the input $u$ is scalar. Triangle block
        symbols are multiplications of the signal by the factor in the block.}}
  \label{fig:block}
\end{figure}

Time-delayed feedback control was originally proposed by Pyragas in
1992 as a tool for discovery of unstable periodic orbits (one frequent
building block in nonlinear systems with chaotic dynamics or multiple
attractors) in experimental nonlinear dynamical systems
\cite{P92}. Pyragas proposed that one take the output $x(t)\in\R^n$ of
a dynamical system and feed back in real time the difference between
this output and the output time $T$ ago into an input
$u(t)\in\R^{n_u}$ of the system (multiplied by some \emph{control
  gains} $K^T\in\R^{n_u\times n}$):
\begin{align}\label{eq:intro:tdf}
  u(t)=K^T[x(t-T)-x(t)]\mbox{.}
\end{align}
\noindent
In a first experimental demonstration, Pyragas and
Tama\v{s}evi\v{c}ius successfully identified and stabilised an
unstable periodic orbit in a chaotic electrical circuit
\cite{PT93a}. Socolar \emph{ et al} in 1994 \cite{SSG94} introduced a
generalisation of time-delayed feedback (which is often used in place
of \eqref{eq:intro:tdf} and is implemented as shown in
Figure~\ref{fig:block} as a block diagram):
\begin{align}
  \label{eq:intro:etdf}
  \begin{split}
    u(t)&=\ K^T[\tilde x(t)-x(t)]\mbox{,\quad where} \\
    \tilde
    x(t)&=(1-\epsilon) \tilde x(t-T)+\epsilon x(t-T)\mbox{,}
  \end{split}
\end{align}
and $\epsilon\in(0,1]$, called \emph{extended time-delayed
  feedback}. If $\epsilon=1$, feedback law \eqref{eq:intro:etdf}
reduces to time-delayed feedback \eqref{eq:intro:tdf}, if $\epsilon=0$
feedback law \eqref{eq:intro:etdf} \rev{degenerates} to classical linear
feedback with a fixed $T$-periodic reference signal $\tilde x(t)$ (see
below \eqref{eq:intro:classical} for a discussion).  \rev{Note that,
  for example in \cite{SSG94}, the variable $\tilde x(t)$ was
  eliminated in the mathematical discussion by writing
  \begin{displaymath}
    u(t)=K^T\left(\epsilon\left[\sum_{j=1}^\infty (1-\epsilon)^{j-1}x(t-jT)\right]-x(t)\right)\mbox{.}
\end{displaymath}
While this would suggest that knowledge of all history of $x$ is required to
initialise the system, in the experiment the feedback control
was implemented as shown in the block diagram in
Figure~\ref{fig:block}, which is equivalent to \eqref{eq:intro:etdf}.}
By construction of the feedback laws \eqref{eq:intro:tdf} and
\eqref{eq:intro:etdf}, \rev{for $\epsilon>0$} every periodic orbit of
period $T$ of the dynamical system with feedback control is also a
periodic orbit of the uncontrolled system ($u=0$).\footnote{\rev{For
    $\epsilon\neq0$, \eqref{eq:intro:etdf} with $T$-periodic $\tilde
    x$ implies that $\tilde x=x$ for all $t$.}} However, the stability
of the periodic orbit may change from unstable without control to
asymptotically stable with control for appropriately chosen gains $K$.

The delayed terms $x(t-T)$ and $\tilde x(t-T)$ make extended
time-delayed feedback control different from the classical linear
feedback control, which has the form
\begin{align}
  \label{eq:intro:classical}
  u(t)=K^T[x_*(t)-x(t)]\mbox{,}
\end{align}
where $x_*(t)$ is, for example, a known unstable periodic orbit of the
dynamical system governing $x$. While the goal of
\eqref{eq:intro:classical} is to stabilise a known \emph{reference
  output} (in this case a periodic orbit), time-delayed feedback is
able to stabilise and, thus, \emph{find} a-priori unknown periodic
orbits. For this reason time-delayed feedback originated, and has
found most attention, in the physics and science community, rather
than in the control engineering community. It can be used to discover
features of nonlinear dynamical systems inaccessible in conventional
experiments, such as unstable equilibria, periodic orbits and their
bifurcations, non-invasively. A few examples where time-delayed
feedback (or its extended version) have been successfully used are:
control of chemical turbulence \cite{KBPOMBRE01}, all-optical control
of unstable steady states and self-pulsations in semiconductor lasers
\cite{UBBWH04,WSH08,SHWSH06}, control of neural synchrony
\cite{PHT05,PHT06,SHHD09}, control of the Taylor-Couette flow
\cite{LWP01}, atomic force microscopy \cite{YH06} and (with further
modifications) systematic bifurcation analysis in mechanical
experiments in mechanical engineering \cite{SGNWK08,BMB12,BSESTS2014}.

One difficulty for time-delayed feedback is that there are until now
no general statements guaranteeing the existence of stabilising
control gains $K$ under some genericity condition on the dynamical
system governing $x$ and its input $u$, such as controllability. This
is in contrast to the situation for classical linear feedback control
\eqref{eq:intro:classical}, where the following is known \cite{MSA04}:
if the periodic orbit $x_*$ is linearly controllable by input $u$ in
$p$ periods (this is a genericity condition) then one can assign its
period-$pT$ monodromy matrix to any matrix with positive determinant
by $pT$-periodic feedback gains $K(t)^T\in\R^{n_u\times n}$.

The greater level of difficulty for (extended) time-delayed feedback
is unsurprising since the feedback-controlled system acquires memory.
Let us assume that the measured quantity $x$ is governed by an
ordinary differential equation (ODE) $\dot x(t)=f(x(t),u(t))$ (which
is autonomous without control ($u=0$) and non-autonomous with
classical feedback control \eqref{eq:intro:classical}). \rev{Then $x$
  and $\tilde x$ will be governed by a delay differential equation
  (DDE) if $u$ is given by time-delayed feedback \eqref{eq:intro:tdf},
  or by an ODE coupled to a difference equation if $u$ is given by
  extended time-delayed feedback \eqref{eq:intro:etdf} with
  $\epsilon\in(0,1)$ (we will refer to both cases simply as
  DDEs). This means that the initial value for both, $x$ and $\tilde
  x$, is a history segment, a function on $[-T,0]$ with values in
  $\R^n$.} In DDEs periodic orbits have infinitely many Floquet
multipliers.\footnote{Floquet multipliers are the eigenvalues of the
  linearisation of the time-$T$ map along the periodic
  orbit.} Section~\ref{sec:review} will review the development of
analysis for the time-delayed feedback laws \eqref{eq:intro:tdf} and
\eqref{eq:intro:etdf}. This paper proves a first simple generic
stabilisability result for extended time-delayed feedback control
\eqref{eq:intro:etdf} with time-periodic gains $K(t)$ (similar to
results for classical linear feedback control).
\paragraph{Main result} 
The following theorem states that the classical approach to periodic
feedback gain design by Brunovsky \cite{B69a} can be applied to make
\eqref{eq:intro:etdf} stable in the limit of small $\epsilon>0$ in the
simplest and most common case of a scalar input $u$ (thus, $n_u=1$)
and linear controllability of the periodic orbit by an input at a
single time instant.
\begin{theorem}[Generic stabilisability with extended time-delayed feedback]
  \label{thm:main}
  Assume that the dynamical system
  \begin{align}
    \label{eq:intro:odeu}
    \dot x(t)&=f(x(t),u(t))\mbox{\quad \textup{(}$f:\R^n\times \R\mapsto\R^n$ smooth\textup{)}
      }
  \end{align}
  with $u=0$ has a periodic orbit $x_*(t)$ of period $T>0$, and assume
  that the monodromy matrix\footnote{The monodromy matrix $P_0$
    is the solution $y$ at time $T$ of the linear differential
    equation $\dot y(t)=\partial_xf(x_*(t),0)y(t)$ with initial value
    $y(0)=\id$ ($\id$ is the identity matrix).} $P_0$ of $x_*$ from time $0$ to $T$ is controllable with
  $b_0=\partial_uf(x_*(0),0)$ \textup{(}that is,
  $\det[b_0,P_0b_0,\ldots,P_0^{n-1}b_0]\neq0$\textup{)}.

  Then there exist gains $K_0\in\R^n$ such that $x_*$ as a periodic
  orbit of the feedback controlled system \eqref{eq:intro:odeu} with
  \textup{(}see below for the definition of the function $\Delta_\delta$\textup{)}
  \begin{align}
    \label{eq:intro:etdf:u}
    u(t)&=\Delta_\delta(t)K_0^T[\tilde x(t)-x(t)]\mbox{,} &
    \tilde x(t)&=(1-\epsilon)\tilde x(t-T)+\epsilon x(t-T)
  \end{align}
  has one simple Floquet multiplier at $1$ and all other Floquet
  multipliers inside the unit circle for all sufficiently small
  $\epsilon$ and $\delta$.
\end{theorem}
The function $\Delta_\delta$ is zero except for a short interval of
length $\delta$ every period $T$ such that the feedback $u$ has the
form of a short but large near-impulse:\footnote{The notation
  $t\vert_{\mod[0,T)}$ refers to the number $\tau\in[0,T)$ such that
  $(t-\tau)/T$ is an integer.}
\begin{align}
  \label{eq:intro:delta}
  \Delta_\delta(t)&=
  \begin{cases}
    1/\delta & \mbox{if $t\vert_{\mod[0,T)}\in[0,\delta]$}\\
    0 & \mbox{if
      $t\vert_{\mod[0,T)}\notin[0,\delta]$.}
  \end{cases}
\end{align}
\begin{figure}[ht]
  \centering
  \includegraphics[scale=0.5]{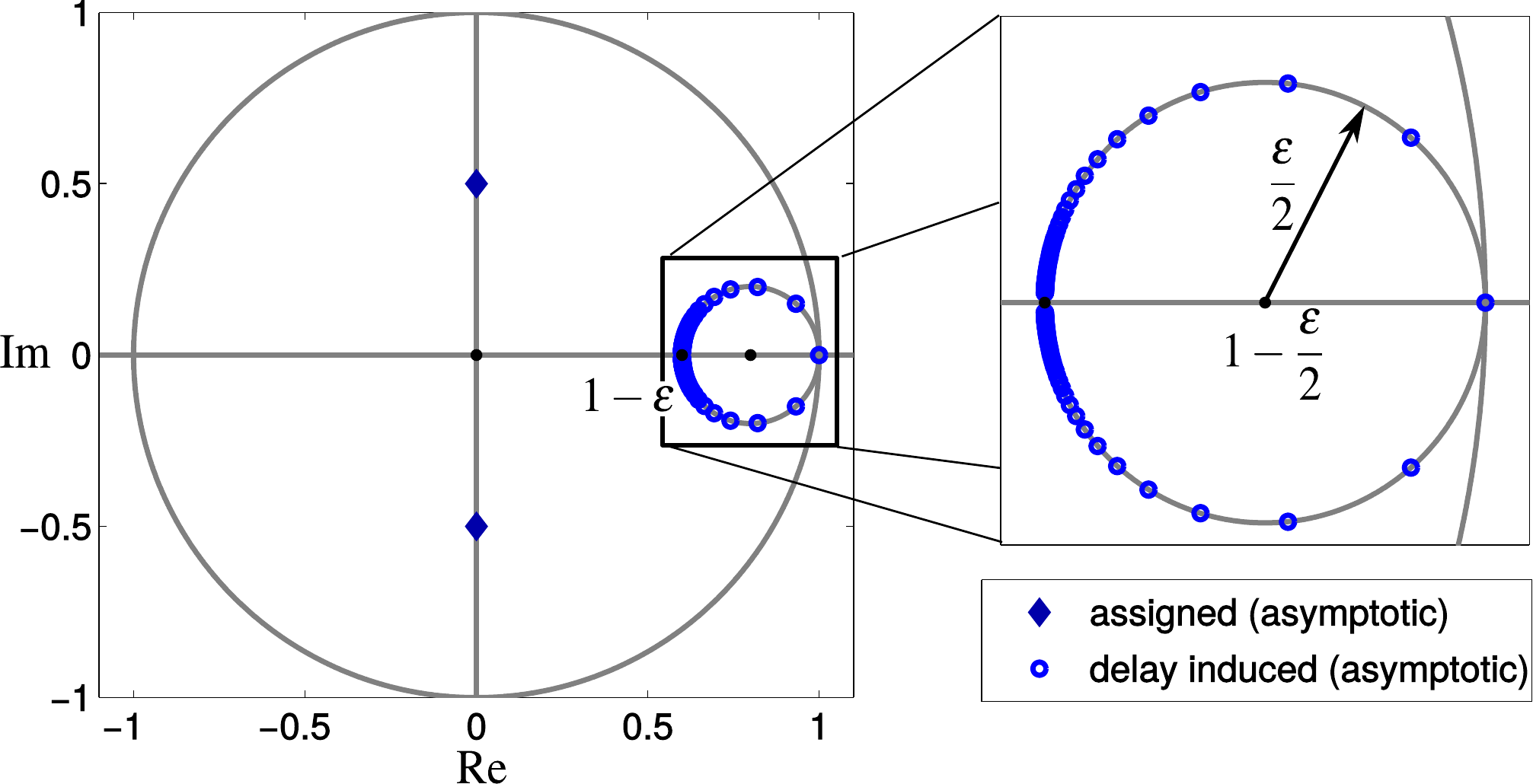}
  \caption{Illustration of Floquet multiplier spectrum for extended
    time-delayed feedback with single input. Using the appropriate
    control gains $K_0$, $n$ Floquet multipliers can be freely
    assigned up to determinant restrictions ($n=2$ in the illustrated
    case). The other Floquet multipliers lie on a circle of radius
    $\epsilon/2$ around $1-\epsilon/2$, accumulating at
    $1-\epsilon$. This spectrum is achieved asymptotically for
    sufficiently short and strong impulses ($\delta\ll1$) and small
    $\epsilon$. A simple trivial multiplier at $1$ is always present.}
  \label{fig:illu}
\end{figure}
\paragraph{Remarks --- Constant gains}
The gains as constructed are periodic. This is to be expected since
there are no general results for constant gains $K\in\R^n$ for the
classical linear feedback case \eqref{eq:intro:classical},
either. Furthermore, simple examples show that the above statement can
definitely not be made when we restrict ourselves to constant gains in
\eqref{eq:intro:etdf:u}: $u(t)=K^T[\tilde x(t)-x(t)]$ with
$K\in\R^n$. See Section~\ref{sec:hopf} for an example.

\paragraph{Properties of the spectrum of the linearisation} \rev{ (see
  also Figure~\ref{fig:illu} for an illustration) The claim of
  Theorem~\ref{thm:main} is about linear stability of the periodic
  orbit $(x(t),\tilde x(t))=(x_*(t),x_*(t))$ of
  \eqref{eq:intro:odeu}--\eqref{eq:intro:etdf:u}. Thus, we have to
  consider the problem \eqref{eq:intro:odeu}--\eqref{eq:intro:etdf:u},
  linearised in $(x(t),\tilde x(t))=(x_*(t),x_*(t))$:
  \begin{equation}\label{eq:intro:lin}
    \begin{split}
      \dot x(t)&=A(t)x(t)+b(t)\Delta_\delta(t)K_0^T[\tilde
      x(t)-x(t)]\mbox{,}  \\
      \tilde x(t)&=(1-\epsilon)\tilde x(t-T)+\epsilon x(t-T)\mbox{,}
    \end{split}
  \end{equation}
  where $ A(t)=\partial_xf(x_*(t),0)$ and
  $b(t)=\partial_uf(x_*(t),0)$. 

  The gains $K_0$ are identical to those chosen by Brunovsky
  \cite{B69a} for the classical feedback spectrum assignment problem
  (note that Brunovsky made weaker assumptions on $A(t)$ and $b(t)$
  than Theorem~\ref{thm:main}). One can choose the gains $K_0$ to
  place the $n$ Floquet multipliers $\lambda_k$ ($k=1,\ldots,n$) of
  \begin{displaymath}
    \dot x=A(t)x-b(t)\Delta_\delta(t)K_0^Tx\mbox{,}
  \end{displaymath} 
  anywhere inside the unit circle subject to the restriction that
  they have to be eigenvalues of a real matrix with positive
  determinant.

  However, DDEs such as \eqref{eq:intro:lin} may have infinitely many
  Floquet multipliers. Theorem~\ref{thm:main} rests on a perturbation
  argument for small $\epsilon>0$ for the other, delay-induced,
  Floquet multipliers.\footnote{The perturbation is \emph{not} a
    small-delay perturbation since the delay $T$ and one coefficient
    in front of the delay, $1-\epsilon$, are not small.} At
  $\epsilon=0$ the difference equation for $\tilde x$ in
  \eqref{eq:intro:lin} simplifies to $\tilde x(t)=\tilde
  x(t-T)$. Thus, an arbitrary initial history $\tilde x$ with period
  $T$ will not change under the time-$T$ map of
  \eqref{eq:intro:lin}. This results for $\epsilon=0$ in a spectrum of
  \eqref{eq:intro:lin} consisting of
  \begin{itemize}
  \item the finitely many assigned Floquet multipliers $\lambda_k$
    ($k=1,\ldots,n$) as determined by the gains $K_0$, and (assuming
    all $\lambda_k\neq1$)
  \item the spectral point $\lambda_\infty=1$ with an infinite-dimensional
    eigenspace. Specifically, if we choose the space of continuous
    functions $C([-T,0];\R^n\times\R^n)$ as phase space for
    \eqref{eq:intro:lin} then, for
    $\epsilon=0$, the eigenspace for $\lambda_\infty=1$ is
    \begin{multline*}
      \Bigl\{(x,\tilde x)\in C([-T,0];\R^n\times\R^n): \quad
      x(0)=x(-T),\quad
      \tilde x(0)=\tilde x(-T), \\
      \dot x=A(t)x+b(t)\Delta_\delta(t)K_0^T\left[\tilde
        x(t)-x(t)\right]\Bigr\}\mbox{.}
    \end{multline*}
    Note that, since $\lambda_k\neq1$ for $k=1,\ldots,n$, the ODE
    $\dot x=A(t)x+b(t)\Delta_\delta(t)K_0^T\left[\tilde
      x(t)-x(t)\right]$ has a unique periodic solution $x$ for all
    periodic functions $\tilde x$. This means that for every
    $T$-periodic $\tilde x$ there is an eigenvector for
    $\lambda_\infty=1$ with this $\tilde x$-component.
  \end{itemize}
  The general theory for DDEs \cite{HL93} ensures that for positive
  (small) $\epsilon$ the Floquet multipliers $\lambda_k$
  ($k=1,\ldots,n$) are only slightly perturbed, and that the
  infinitely many Floquet multipliers emerging from $\lambda_\infty$
  accumulate to the spectrum of the \emph{essential part}, the
  difference equation in \eqref{eq:intro:lin} with the $\tilde x$
  terms only: $\tilde x(t)=(1-\epsilon) \tilde x(t-T)$. Specifically,
  the only accumulation point of the spectrum of \eqref{eq:intro:lin}
  for $\epsilon\in(0,1)$ is at $1-\epsilon$ and the stability of
  \eqref{eq:intro:lin} is determined by the location of the Floquet
  mulitpliers emerging from the perturbation of $\lambda_\infty$ (of
  which at most finitely many can lie outside the unit circle). The
  detailed analysis in Section~\ref{sec:lin:etdf} will show that for
  small $\epsilon>0$ the Floquet multipliers emerging from
  $\lambda_\infty$ lie close to a circle of radius $\epsilon/2$ around
  $1-\epsilon/2$, inside the unit circle (except for the unit Floquet
  multiplier), as shown in Figure~\ref{fig:illu} for $n=2$.}

\paragraph{Trivial multiplier} The eigenvector to the trivial
multiplier $1$ is $\dot x_*(0)$, corresponding to the linearised phase
shift (for every $s\in\R$, $t\mapsto x_*(t+s)$ is also a solution of
the system with extended time-delayed feedback
\eqref{eq:intro:odeu},\,\eqref{eq:intro:etdf:u}). Section~\ref{sec:auto}
gives a modification of Theorem~\ref{thm:main} with a function
$\Delta_\delta$ depending on $x(t)$ instead of $t$ switching the gains
on and off. Then the feedback controlled system becomes autonomous. In
this modified system with autonomous (but nonlinear) extended
time-delayed feedback the periodic orbit $x_*$ is asymptotically
stable in the classical sense.

\paragraph{Timing of impulse} In \eqref{eq:intro:delta} we chose the
timing of the impulse (the part of the period $[0,T)$ where
$\Delta_\delta$ is non-zero) as $[0,\delta]$ without loss of
generality. The genericity condition in its most general form requires
that there must be a time $t\in[0,T]$ such that the monodromy matrix
from $t$ to $t+T$ and $\partial_uf(x_*(t),0)$ are controllable. As the
uncontrolled system is autonomous, we can shift the phase of the
periodic orbit $x_*$, considering $x_*(t+\cdot)$ instead of $x_*$.

\paragraph{Practical considerations} The result gives precise control
over the Floquet multipliers in the limit of small $\delta$ and
$\epsilon$. For small $\delta$ the feedback control corresponds to a
sharp kick once per period, which is not practical for strongly
unstable periodic orbits. However, the gains found with the help of
Theorem~\ref{thm:main} provide a feasible starting point for
optimisation-based spectrum assignment methods (\emph{continuous pole
  placement}) as constructed by Michiels \emph{et al}
\cite{MEVR02,MN14} and adapted to time-delayed
feedback~\eqref{eq:intro:tdf} \cite{LHFGFS11,PP13,PP14}. In the
context of continuation one can combine the gains provided by
Theorem~\ref{thm:main} as starting points, continuous pole placement,
and the automatic adjustment of the time delay $T$ demonstrated in
\cite{PP11,NP12} to create a feedback control that non-invasively
tracks a family of periodic orbits in a system parameter.

\section{Review: analysis of (extended) time-delayed
  feedback}
\label{sec:review}
The initial proposals of time-delayed feedback \eqref{eq:intro:tdf} and
its extended version \eqref{eq:intro:etdf} were accompanied with
demonstrations in simulations and experiments, showing that this type
of feedback control can be successful \cite{P92,PT93a,SSG94}, but not
with general necessary or sufficient conditions for applicability or
with constructive ways to design the feedback gains.

However, it was quickly recognised that time-delayed feedback can be
applied to periodic orbits that are weakly unstable due to a period
doubling bifurcation or torus bifurcation \cite{JBORB97,SS07}. Hence,
time-delayed feedback is often associated with control of chaos,
because it can be used to suppress period doubling cascades. However,
general sufficient criteria were rather restrictive \cite{N04},
requiring full access to the state ($x$ governed by $\dot
x(t)=f(x(t))+u(t)$ with $u\in\R^n$). A first general result was
negative, the so-called \emph{odd number limitation} for periodically
forced systems \cite{NU98}, showing that extended time-delayed feedback
cannot stabilise periodic orbits in periodically forced systems with
an odd number of Floquet multipliers $\lambda$ with $\Re\lambda>1$
(and no Floquet multiplier at $1$). This theoretical limitation is
not a severe restriction in practice since one can extend the
uncontrolled system with an artificial unstable degree of freedom
before applying time-delayed feedback \cite{P01}. Fiedler \emph{et al}
showed that this limitation does not apply to autonomous periodic
orbits \cite{FFGHS07,FS11}. Since then general results have been
proven for weakly unstable periodic orbits with a Floquet multiplier
close to $1$ (but larger than $1$, \cite{PN13}), or near subcritical Hopf
bifurcations \cite{BPS11,PBS13}. A review of developments up to 2010
is given in \cite{SHFD10}. 

An extension of the odd number limitation to autonomous periodic
orbits (with trivial Floquet multiplier) was given by Hooton \& Amann
\cite{HA12,AH13} for both, time-delayed feedback \eqref{eq:intro:tdf}
and its extension \eqref{eq:intro:etdf}. However, these limitations
merely impose restrictions on the gains $K$. They do not rule out
feedback stabilisability a priori (which is in contrast to the
statements about periodic orbits in forced systems).

\section{Spectrum of linearisation for extended time-delayed
  feedback-controlled system}
\label{sec:lin:etdf}
Let us consider a feedback controlled dynamical system with extended
time-delayed feedback control and arbitrary time-dependent gains
$K(t)\in\R^n$:
\begin{align}
  \dot x(t)&=f(x(t),u(t))\mbox{,}\label{eq:etdf:x}\\
  \label{eq:etdf:u}
    u(t)&=K(t)^T[\tilde x(t)-x(t)]\mbox{,}\\
    \tilde x(t)&=(1-\epsilon)\tilde x(t-T)+\epsilon x(t-T)\mbox{.}
  \label{eq:etdf:xtilde}
\end{align}
This system is governed by an ordinary differential equation (ODE)
without control ($u=0$) and a delay differential equation (DDE) with
control.  We assume that the uncontrolled system $\dot x(t)=f(x(t),0)$
has a periodic orbit $x_*$ of period $T$. This periodic orbit $x_*$ is
also a periodic orbit of \eqref{eq:etdf:x}--\eqref{eq:etdf:xtilde} if
$\epsilon>0$: $x(t)=\tilde
x(t)=x_*(t)$. System~\eqref{eq:etdf:x}--\eqref{eq:etdf:xtilde} is a
DDE with the phase space
\begin{align*}
  \left\{(x,\tilde x)\in C([-T,0];\R^n\times\R^n):\tilde
  x(0)=(1-\epsilon)\tilde x(-T)+\epsilon x(-T)\right\}\mbox{.}
\end{align*}
\rev{Hale \& Verduyn-Lunel \cite{HL93} treated DDEs of the type of
  system~\eqref{eq:etdf:x}--\eqref{eq:etdf:xtilde} (which contains
  difference equations) as part of their discussion of neutral DDEs.}
The essential part of the semiflow generated by
\eqref{eq:etdf:x}--\eqref{eq:etdf:xtilde} is governed by the part of
\eqref{eq:etdf:xtilde} containing $\tilde x$: $\tilde
x(t)=(1-\epsilon)\tilde x(t-T)$, which is linear and has spectral
radius $1-\epsilon$. Thus, it fits into the scope of the theory as
described in the textbook by Hale \& Verduyn-Lunel
\cite{HL93}. Specifically, the asymptotic stability of the periodic
orbit given by $x(t)=\tilde x(t)=x_*(t)$ is determined by the point
spectrum of the linearisation of
\eqref{eq:etdf:x}--\eqref{eq:etdf:xtilde}.  Hence, the periodic orbit
$x_*$ is stable if all Floquet multipliers of the linearisation along
$x_*$ except the trivial multiplier $1$ are inside the unit circle
(and the trivial Floquet multiplier $1$ is simple). We denote the
monodromy matrix\footnote{Thus, $P(\mu)$ is defined as the solution
  $y$ at time $T$ of the linear differential equation $\dot
  y(t)=[A(t)-\mu b(t)K(t)^T]y(t)$ with initial value $y(0)=\id$ ($\id$
  is the identity matrix).} of
\begin{align}\label{eq:monodromy}
  \dot x&=[A(t)-\mu b(t)K(t)^T]x(t)\mbox{,}&\mbox{where\ }
  A(t)&=\partial_xf(x_*(t),0)\mbox{,}& b(t)&=\partial_uf(x_*(t),0)
\end{align}
for $\mu\in\C$ by $P(\mu)$. Thus, the monodromy matrix of the
uncontrolled system $\dot x(t)=A(t)x(t)$ equals $P(0)$, which we denote by
\begin{align}
  \label{eq:p0def}
  P_0=P(0)\mbox{.}
\end{align}
With this definition of $P(\mu)$, Floquet multipliers of the
linearisation of \eqref{eq:etdf:x}--\eqref{eq:etdf:xtilde} in $x_*$
different from $1-\epsilon$ are given as roots of
\begin{align*}
  h(\lambda;\epsilon):=
  \det\left[\lambda\id-P\left(1-
      \frac{\epsilon}{\lambda-(1-\epsilon)}\right)\right]
\end{align*}
($\id$ is the identity matrix; \rev{see
  Section~\ref{sec:proof:floquet} for detailed proof}).  The following
lemma states that the gains $K(t)$ can only stabilise a periodic orbit
$x_*$ with extended time-delayed feedback and small $\epsilon$, if
they are stabilising with classical linear feedback (that is, when
replacing the recursively determined signal $\tilde x$ by the target
orbit $x_*$: $u(t)=K(t)[x_*(t)-x(t)]$). (Recall that
$A(t)=\partial_xf(x_*(t),0)$, $b(t)=\partial_uf(x_*(t),0)$.)
\begin{lemma}[Extended time-delayed feedback stabilisation implies
  classical stabilisation]
  \label{thm:etdf2classical}
  If the linear system
  \begin{align}\label{eq:classical:feedback}
    \dot x(t)=[A(t)-b(t)K(t)^T]x(t)
  \end{align}
  has at least one Floquet multiplier outside the unit circle, then
  there exists a $\epsilon_{\max}\in(0,1)$ such that the periodic
  orbit $x_*$ is unstable for the extended time-delayed feedback
  \eqref{eq:etdf:x}--\eqref{eq:etdf:xtilde} for all
  $\epsilon\in(0,\epsilon_{\max})$.
\end{lemma}
\paragraph{Proof} The Floquet multipliers of
\eqref{eq:classical:feedback} are given as roots of
$h(\lambda;0)=\det(\lambda\id-P(1))$. We denote the root with modulus
greater than $1$ by $\lambda_0$ such that
$h(\lambda_0;0)=0$. Consequently, for all $\lambda$ in the ball
$B_r(\lambda_0)$, where $r=(|\lambda_0|-1)/2$, the difference
$h(\lambda;\epsilon)-h(\lambda;0)$ is uniformly bounded and analytic
for all $\epsilon\in(0,1)$ and all $\lambda$ in
$B_r(\lambda_0)$. Since $\lambda_0$ must have finite multiplicity as a
root of $h(\cdot;0)$, $h(\lambda;\epsilon)$ must have a root in
$B_r(\lambda_0)$ for sufficiently small $\epsilon>0$ (say,
$\epsilon\in(0,\epsilon_{\max})$), too. By choice of $r$ this root
lies outside of the unit circle. \eop{Lemma~\ref{thm:etdf2classical}}

Lemma~\ref{thm:etdf2classical} shows that gains $K(t)$ that stabilise
with extended time-delayed feedback with small $\epsilon$ also have to
feedback-stabilise in the classical sense. Since $K$ is an arbitrary
periodic function there are many ways to construct gains for the
classical linear feedback control $u(t)=K(t)[x_*(t)-x(t)]$ for
periodic orbits $x_*$ \cite{MSA04}. We choose the approach comprehensively treated
by Brunovsky \cite{B69a}, which is particularly amenable to analysis
in the extended time-delayed feedback case and for which one can then
prove the converse of Lemma~\ref{thm:etdf2classical}:
\begin{quote}
  \emph{extended time-delayed feedback is stabilising for the same
    gains for which Brunovksy's approach is stabilising  the classical
    feedback control.}
\end{quote}
\paragraph{Near-impulse feedback and its parametrised monodromy matrix}
We pick state feedback control in the form of a single large
but short impulse. That is, we consider a short time $\delta\in(0,T)$
and define the linear feedback control
\begin{align}
  \label{eq:genfeedback}
  u_\delta(t;y)&=\Delta_\delta(t)K_0^Ty\mbox{, where} &
  \Delta_\delta(t)&=
  \begin{cases}
    1/\delta & \mbox{if $t\vert_{\mod[0,T)}\in[0,\delta]$}\\
    0 & \mbox{if $t\vert_{\mod[0,T)}\notin[0,\delta]$.}
  \end{cases}
\end{align}
where 
$t\vert_{\mod[0,T)}$ is the number $\tau\in[0,T)$ such that
$(t-\tau)/T$ is an integer, and $K_0\in\R^n$ is a vector of constant
control gains. Let us first look at classical feedback
$u(t)=\Delta_\delta(t)K_0^T[x_*(t)-x(t)]$ (where we assume that we know the periodic
orbit $x_*$). Using feedback law \eqref{eq:genfeedback} the
feedback controlled system reads
\begin{equation}
  \label{eq:lincont}
  \dot x(t)=f\left(x(t),\Delta_\delta(t)K_0^T[x_*(t)-x(t)]\right)\mbox{.}
\end{equation}
We define the nonlinear time-$T$ map $X(x;\delta,K_0)$ as the solution
at time $T$ (the period of the periodic orbit $x_*$) of
\eqref{eq:lincont} when starting from $x$ at time $0$ (including the
dependence on parameters $\delta$ and $K_0$ as additional arguments of
$X$).  Then, for small deviations $y_0$ from $x_*(0)$, the map
$X(\cdot;\delta,K_0)$ has the form
$X(x_*(0)+y_0;\delta,K_0)=y(T)+O(\|y_0\|^2)$, where $y$ satisfies the
linear differential equation (recall that
$A(t)=\partial_xf(x_*(t),0)$, $b(t)=\partial_uf(x_*(t),0)$)
\begin{align}
  \label{eq:tlin}
  \dot y(t)&=\left[A(t)-b(t)\Delta_\delta(t)K_0^T\right]y(t)\mbox{,}& 
  y(0)&=y_0\mbox{,}
\end{align}
and the term $O(\|y_0\|^2)$ is uniformly small (including its
derivatives) for all $\delta$ . Let us introduce a complex parameter
$\mu$ into \eqref{eq:tlin}, which will become useful later in our
consideration of extended time-delayed feedback: define for a general
complex $\mu$ with $|\mu|\leq C$ (with an arbitrary fixed
$C>0$) the linear ODE
\begin{align}
  \label{eq:mulin}
  \dot y(t)&=\left[A(t)-\mu b(t)\Delta_\delta(t)K_0^T\right]y(t)\mbox{,}& 
  y(0)&=y_0\mbox{.}
\end{align}
Denote the monodromy matrix of \eqref{eq:mulin} from $t=0$ to $t=T$ by
$P(\mu;\delta,K_0)$ to keep track of its dependence on the parameters
$\delta\in(0,T)$ and $K_0\in\R^n$. Thus, $P(\mu;\delta,K_0)$ refers to
the same monodromy matrix as $P(\mu)$, defined by
\eqref{eq:monodromy}, for the special case
$K(t)=\Delta_\delta(t)K_0$. Then $P(\mu;\delta,K_0)$ satisfies
\begin{align}\label{eq:Plimit}
  P(\mu;\delta,K_0)=
  P_0\exp(-b(0)K_0^T\mu)+O(\delta)\mbox{,}
\end{align}
where the error term $O(\delta)$ is uniform for $|\mu|\leq C$ and
bounded $\|K_0\|$, including its derivatives with respect to all
arguments.  Hence, we can extend the definition of
$P(\mu;\delta,K_0)$ to $\delta=0$:
\begin{align}
  P(\mu;0,K_0)=&\ \lim_{\delta\to0}P(\mu;\delta,K_0)=
  P_0\exp(-b(0)K_0^T\mu)\mbox{,}    \nonumber\\
  =&\ P_0[\id-\sigma(\mu)b(0)K_0^T]\mbox{\quad where}\label{eq:p0deltadef}\\
  \sigma(\mu)=&\   
  \begin{cases}
    \cfrac{\exp\left(\mu K_0^Tb(0)\right)-1}{K_0^Tb(0)} 
    &\mbox{if $K_0^Tb(0)\neq0$}\\
    \mu &\mbox{if $K_0^Tb(0)=0$.}
  \end{cases}\label{eq:sigmadef}
\end{align}
The limit is uniform for all $\mu$ with modulus less than $C$. For
$\mu=0$, $P$ is the monodromy matrix $P_0$ of the uncontrolled system,
and, thus, independent of $\delta$ and $K_0$.

\paragraph{Approximate spectrum assignment for finitely many Floquet
  multipliers}
The control \eqref{eq:genfeedback} is a simplification of the general
case of finitely many (at most $n$) short impulses treated in
\cite{B69a}. Feedback of type \eqref{eq:genfeedback} permits us to
assign arbitrary spectrum approximately under the assumption that the
pair $(P_0,b(0))$ is controllable (recall that, according to the
definition of $P_0$ in \eqref{eq:p0def}, $P_0$ is the monodromy matrix
of the uncontrolled system $\dot x(t)=f(x(t),0)$ along the periodic
orbit $x_*$). This is a stronger assumption than the assumption made
in \cite{B69a}, but it is still a genericity assumption.
\begin{lemma}[Approximate spectrum assignment for classical state
  feedback control, simplified from \cite{B69a}]\label{thm:stab}
  Let $r>0$ be arbitrary.
  If the pair $(P_0,b(0))$ is controllable \textup{(}that is, the $n\times
  n$ controllability matrix
  \begin{math}
    \left[ b(0),P_0b(0),\ldots, P_0^{n-1}b(0)\right]
  \end{math}
  is regular\textup{)}, then there exist a $\delta_{\max}>0$ and a vector of
  control gains $K_0\in\R^n$ in \eqref{eq:genfeedback} such that all
  Floquet multipliers of $x_*$ for the differential equation
  \eqref{eq:lincont} have modulus less than $r$ for all $\delta\in
  (0,\delta_{\max})$, where $\Delta_\delta$ is as defined in
  \eqref{eq:genfeedback}.
\end{lemma}
Note that the vector $K_0$ can be chosen independent of the
$\delta\in(0,\delta_{\max})$, but it may depend on the radius $r$ into
which one wants to assign the spectrum. This result follows from
classical linear feedback control theory (\cite{B69a} proves a more
general result). 
In short, linear feedback control theory \cite{B69a} makes the
following argument (thus, proving Lemma~\ref{thm:stab}): the
linearisation of $X$ with respect to its initial condition can be
expanded in $\delta$ as
\begin{align*}
  \partial_xX(x_*(0);\delta,K_0)=P(1;\delta,K_0)=
  P_0\exp(-b(0)K_0^T)+O(\delta)
\end{align*}
(where $P(\cdot;\delta,K_0)$ was the generalised monodromy matrix
defined for \eqref{eq:mulin}).  Since $\det P_0$ is positive we can
for every matrix $R$ with positive determinant find a vector $K_0$
such that $\spec R=\spec(P_0\exp(-b(0)K_0^T))$ (using the assumption
of controllability; see auxiliary Lemma~\ref{thm:specass}, which is a
special case from the more general treatment in \cite{B69a}, and
\cite{S16supp} for a Matlab implementation). Hence, if we choose the
spectrum of $R$ inside a circle $B_{r/2}(0)$ of radius $r/2$ around
$0$, then the spectrum of $\partial_x X(x_*(0);\delta,K_0)$ is also
inside $B_r(0)$ for sufficiently small $\delta>0$.

\paragraph{Approximate spectrum for extended time-delayed feedback}
We fix the control gains $K_0$ such that $\lim_{\delta\to
  0}P(1;\delta,K_0)=P_0\exp(-b(0)K_0^T)$ has all eigenvalues inside
$B_r(0)$ for some $r\in(0,1)$. Consider now again the extended
time-delayed feedback control
\eqref{eq:etdf:x}--\eqref{eq:etdf:xtilde} with the particular choice
of short impulse linear feedback law \eqref{eq:genfeedback}:
\begin{align}
  \label{eq:etdfode}
  \dot x(t)&=f(x(t),\Delta_\delta(t)K_0^T[\tilde x(t)-x(t)])\\
  \tilde x(t)&=(1-\epsilon)\tilde x(t-T)+\epsilon x(t-T)\mbox{,}
  \label{eq:etdfref}
\end{align}
where $\epsilon\in(0,1)$.  

\begin{lemma}[Floquet multipliers of extended time-delayed feedback]
  \label{thm:etdf:stable} 
  Assume that the matrix $P_0\exp(-b(0)K_0^T)$ has all eigenvalues
  inside the ball $B_r(0)$ with $r<1$. Then, for all sufficiently 
  small $\epsilon$ and $\delta$, the periodic orbit $x(t)=\tilde
  x(t)=x_*(t)$ of system \eqref{eq:etdfode},\,\eqref{eq:etdfref} has a
  simple Floquet multiplier $\lambda=1$ and all its other Floquet
  multipliers are inside the unit circle.
\end{lemma}

\paragraph{Outline of proof (details are given in Section~\ref{sec:proof:etdf:stable})}
Eigenvalues $\lambda$ of the linearisation of
\eqref{eq:etdfode}--\eqref{eq:etdfref} are roots of
the function
\begin{align}\label{eq:hdef}
  h(\lambda;\epsilon,\delta)=\det \left[\lambda\id-P
    \left(1-\frac{\epsilon}{\lambda-(1-\epsilon)};\delta,K_0\right)\right]\mbox{.}
\end{align}
Roots of $h$ with a non-small distance from $1-\epsilon$ are close to
the roots of $\det(\lambda\id-P(1;\delta,K_0))$, which are inside the
unit circle by assumption. Roots $\lambda$ of $h$ close to
$1-\epsilon$ with modulus greater than $1-\epsilon/2$ have the form
$\lambda=1-\epsilon+\epsilon/\kappa$ where $|\kappa|$ is bounded away
from $0$ and infinity. The roots $\kappa$ of
$h(1-\epsilon+\epsilon/\kappa;\delta,\epsilon)$ are small
perturbations of the roots $\kappa_{\ell,0}$ of
$\det(\id-P_0-P_0b(0)K_0^T\sigma(\kappa-1))$, where $\sigma$ is as
defined in \eqref{eq:sigmadef}. These roots $\kappa_{\ell,0}$ have the
form
\begin{align}\label{eq:kappa_l0}
  \kappa_{\ell,0}&=1+\frac{2\pi\i\ell}{K_0^Tb(0)}
\end{align}
(if $K_0^Tb(0)\neq0$, otherwise, only a single root $\kappa_{0,0}=1$
exists). The roots $\kappa_{\ell,0}$ have all modulus greather than
unity (except for $\ell=0$, which corresponds to the trivial
eigenvalue $\lambda=1$) such that the corresponding roots
$\lambda_\ell$ of $h$ have modulus smaller than unity.
\eop{Lemma~\ref{thm:etdf:stable}}

\paragraph{Remark --- two types of Floquet multipliers} The proof of
Lemma~\ref{thm:etdf:stable} shows that there are two distinct types of
roots: those approximating the spectrum assigned by the choice of
control gains $K_0$, and those close to $1-\epsilon$ (called
$\lambda_\ell$ above). The roots $\lambda_\ell$ lie close to the
circle of radius $\epsilon/2$ around the center $1-\epsilon/2$ in the
complex plane and have the form
\begin{align*}
  \lambda_\ell&\approx 1-\frac{\epsilon}{2}+
  \frac{\epsilon}{2}
  \left[\frac{K_0^Tb(0)-2\pi\i\ell}{K_0^Tb(0)+2\pi\i\ell}\right]&
  &&\mbox{($\ell\in\Z$).}
\end{align*}
For $\ell=0$, the expression is exact (giving the simple root at
unity), for the others the approximation is sufficiently accurate for
small $\delta$ and $\epsilon$ to ensure that they stay inside the unit
circle. The illustration in Figure~\ref{fig:illu} shows the two
distinct groups for the Hopf normal form example discussed in
Section~\ref{sec:hopf}.

\paragraph{Importance of scalar input and trivial Floquet multiplier}
The proof of Lemma~\ref{thm:etdf:stable} hinges on one argument that
depends on the presence of a trivial Floquet multiplier: we need to
find the roots $s_j$ of $s\mapsto \det(\id-P_0-P_0b(0)K_0^Ts)$ and
then find solutions $\kappa$ of $\sigma(\kappa-1)=s_j$ for all these
roots $s_j$. Since $b(0)K_0^T$ has rank one we know that $s\mapsto
\det(\id-P_0-P_0b(0)K_0^Ts)$ is a first-order polynomial (see
Section~\ref{sec:proof:etdf:stable} for details). The presence of a
trivial Floquet multiplier then ensures that this first-order
polynomial has the root $0$. Hence, $0$ is its only root, restricting
the possible location for the $\kappa_{\ell,0}$ to the list in
\eqref{eq:kappa_l0}. This simple argument would not apply for cases
where the uncontrolled periodic orbit $x_*$ has no trivial Floquet
multiplier, or for control with non-scalar inputs $u$, or for control
with more than one kick per period.

\section{Autonomous feedback control}
\label{sec:auto}
\rev{The feedback control constructed in Lemma~\ref{thm:etdf:stable}
  introduces an explicit time dependence into the system. The
  controlled system has the form
  \begin{equation}\label{eq:etdf:noauto}
    \begin{split}
      \dot x(t)&=f(x(t),\Delta_\delta(t)K_0^T[\tilde x(t)-x(t)])\\
      \tilde x(t)&=(1-\epsilon)\tilde x(t-T)+\epsilon x(t-T)\mbox{,}
    \end{split}
  \end{equation}
  where $\Delta_\delta$ is time-periodic with period $T$, but the
  system still has a Floquet multiplier $\lambda=1$. The neutrally
  stable direction corresponding to this Floquet multiplier is a phase
  shift: if $(x(t),\tilde x(t))=(x_*(t),x_*(t))$ is a periodic orbit
  of \eqref{eq:etdf:noauto} then so is $(x(t),\tilde
  x(t))=(x_*(t+s),x_*(t+s))$ for any $s\in\R$. Hence, the controlled
  system with the gains $K(t)=\Delta_\delta(t)K_0$ is susceptible to
  arbitrarily small time-dependendent perturbations (say, experimental
  disturbances): the phase $s$ of the stabilised solution may drift
  until $\Delta_\delta$ is non-zero at a time $s$ where the gains
  $K_0$ are no longer stabilising. This problem does not occur if,
  instead of applying the feedback $K_0^T[\tilde x(t)-x(t)]$ at a
  fixed time per period, we apply it in a strip in $\R^n$ close to a
  Poincar{\'e} section at $x_*(0)$ (as illustrated in
  Figure~\ref{fig:poincare}), putting a factor depending on $x(t)$ in
  front of $K_0^T[\tilde x(t)-x(t)]$. Specifically, we let the
  function $\Delta_\delta$ not depend explicitly on time $t$ but on a
  function $\tilde t:\R^n\mapsto \R$, where the argument of $\tilde t$
  is $x(t)$.  Then the common notion of asymptotic stability of
  periodic orbits in autonomous dynamical systems applies.
  \begin{figure}[h]
    \centering
    \includegraphics[scale=0.5]{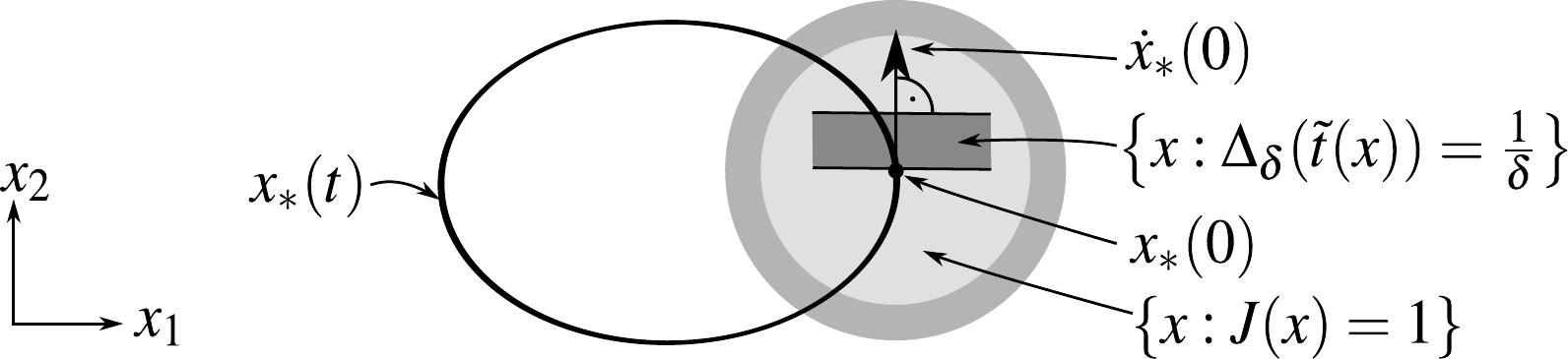}
    \caption{Illustration of choice for strip and Poincar{\'e} section
      where gains should be non-zero.}
    \label{fig:poincare}
  \end{figure}
 One would then always apply
  control near $x_*(0)$ despite phase drift. A possible explicit
  expression for $u$ is
  \begin{align}
    u(t)&=\Delta_{\rho,\delta}(x(t))K_0^T[\tilde
    x(t)-x(t)]\mbox{,}&\mbox{\ where\quad }
    \Delta_{\rho,\delta}(x)&=J_\rho(x)\Delta_\delta(\tilde t(x))\mbox{,}
    \label{eq:auto:u}\\
    J_\rho(x)&=
    \begin{cases}
      1 & \mbox{if $|x-x_*(0)|<\rho$,}\\
      0 & \mbox{if $|x-x_*(0)|>2\rho$,}
    \end{cases}& \tilde t(x)&=\frac{\dot x_*(0)^T}{\dot x_*(0)^T\dot x_*(0)}[x-x_*(0)]\mbox{,}\label{eq:auto:Jt}
  \end{align}
  $\rho>0$ is a small radius and $J_\rho$ is smooth. In
  \eqref{eq:auto:Jt}, $\tilde t(x_*(t))=t+O(t^2)$ for $|t|\ll1$, and
  $J_\rho$ restricts control to the neighborhood of radius $\rho$
  around $x_*(0)$. With $u$ as defined in \eqref{eq:auto:u}, the
  right-hand side of the now autonomous system\begin{align}
  \label{eq:etdfnlin:all}
  \dot x(t)&=f(x(t),\Delta_{\delta,\rho}(x(t))K_0^T[\tilde x(t)-x(t)])\mbox{,} &
  \tilde x(t)&=(1-\epsilon)\tilde x(t-T)+\epsilon x(t-T)
\end{align}
has a right-hand side that depends discontinuously on $x(t)$ (because
$\Delta_\delta$ is discontinuous in its argument. Since the general
mathematical theory for DDEs coupled to difference equations is not
well developed, one may replace the discontinous $\Delta_\delta$ in
\eqref{eq:auto:u} with a smooth approximation of $\Delta_\delta$.
This does not affect the final result, which we can state as a lemma
(see Appendices \ref{sec:ktmod} and \ref{sec:sdK} for the details of
the choice for $\rho$ and the smoothing of $\Delta_{\delta,\rho}$):}
\begin{lemma}[Autonomous stabilisability of periodic orbits with
  extended time-delayed feedback]
  \label{thm:etdf:autostab}
  Assume that the matrix $P_0\exp(-b(0)K_0^T)$, as used in
  Lemma~\ref{thm:etdf:stable}, has all eigenvalues inside the ball
  $B_r(0)$ with $r<1$. Then, for all sufficiently small $\rho$, there
  exist $\epsilon_{\max}>0$ and $\delta_{\max}>0$ such that the
  periodic orbit $x(t)=\tilde x(t)=x_*(t)$ of system
  \eqref{eq:etdfnlin:all} is asymptotically exponentially stable for
  all $\epsilon\in(0,\epsilon_{\max})$ and
  $\delta\in(0,\delta_{\max})$.
\end{lemma}
\paragraph{Remark: other arguments for $\Delta_{\delta,\rho}$}
In \eqref{eq:auto:u} we can replace the argument $x(t)$ of
$\Delta_{\delta,\rho}$ with $\tilde x(t)$, $x(t-T)$ or $\tilde x(t-T)$
without changing the linearisation in $x(t)=\tilde x(t)=x_*(t)$. Thus,
\eqref{eq:etdfnlin:all} successfully stabilises the periodic orbit
$x_*$ also with these modifications.

\paragraph{Robustness}
We assumed perfect knowledge of the periodic orbit $x_*$ and the
right-hand side $f$ in the construction of $K_0$ and
$\Delta_{\delta,\rho}$. However, we know that stable periodic orbits
persist under small perturbations. Thus, for gains near $K_0$ and
functions close to $\Delta_{\delta,\rho}$ the periodic orbit of the
controlled system persists. Due to the non-invasive nature of extended
time-delayed feedback, the periodic orbit of the system with perturbed
$K_0$ and $\Delta_{\delta,\rho}$ is still identical to $x_*$. 

\section{Illustrative example: Hopf normal form}
\label{sec:hopf}
\rev{The construction of gains as described in Section~\ref{sec:auto}
  has been implemented as a Matlab function (publically available at
  \cite{S16supp}, depending on DDE-Biftool
  \cite{ELR02,ELS01,ddebiftoolmanual}). The supplementary material
  demonstrates how one can find stabilising gains for two examples:
  \begin{enumerate}
  \item a family of period-two unstable oscillations around the
    hanging-down position of the parametrically excited pendulum, and
  \item\label{ex:hopf} the unstable periodic orbits in the subcritical
    Hopf normal form.
  \end{enumerate}
  We discuss example~\ref{ex:hopf} in more detail in this section,
  because for this example we can prove that stabilisation with ETDF
  is not possible with constant gains and small $\epsilon$. The
  subcritical Hopf bifurcation has also been used commonly in the
  literature as a benchmark example. Here we choose the Hopf normal
  form with constant speed of rotation (such that in polar coordinates
  the angle $\theta$ satisfies $\dot\theta=1$ and all periodic orbits
  have period $2\pi$). Note that the control constructed by Fiedler
  \emph{et al} \cite{FFGHS07} depended on changing rotation and was
  stabilising only in a small neighborhood of the
  bifurcation. Flunkert \& Sch{\"o}ll \cite{FS11} analysed
  time-delayed feedback control (with $\epsilon=1$) of the subcritical
  Hopf bifurcation completely, but also excluded the case of constant
  rotation and restricted themselves to a small neighbourhood of the
  bifurcation. Thus, even though example~\ref{ex:hopf} is seemingly
  simple, it shows that the method proposed in the paper is able to
  stabilise periodic orbits that are beyond the approaches previously
  suggested in the literature. Without loss of generality we choose a
  linear control input $b=[1,1]^T$ such that the system with control
  has the form:}
\begin{equation}
  \label{eq:hopfu}
  \begin{split}
    \dot x_1 &= p x_1-x_2+x_1[x_1^2+x_2^2]+u\mbox{,}\\
    \dot x_2 &= x_1+px_2+x_2[x_1^2+x_2^2]+u\mbox{,}
  \end{split}
\end{equation}
where $p<0$. This system has for $u=0$ an unstable periodic orbit of
the form $x_*(t)=[r\sin t,-r\cos t]^T$ with radius $r=\sqrt{-p}$ and
period $T=2\pi$. The monodromy matrix $P_0$ for the uncontrolled
system along the periodic orbit $x_*$ equals
\begin{align*}
  P_0&=
  \begin{bmatrix}
    1 & 0 \\ 0 & \exp(-4\pi p)
  \end{bmatrix}\mbox{.}
\end{align*}
Since the derivative of the right-hand side with respect to the
control input equals $b(t)=b=[1,1]^T$, the periodic orbit is
controllable in time $T$. In fact, the pair $(P_0,b)$ is controllable
as required for the applicability of
Lemma~\ref{thm:etdf:stable}. Extended time-delayed feedback control,
applied to a two-dimensional system has the form
\begin{equation}\label{eq:hopf:etdf}
  \begin{split}
    u(t)&=K_1(t)[\tilde x_1(t)-x_1(t)]+K_2(t)[\tilde x_2(t)-x_2(t)]\\
    \tilde x_j(t)&=(1-\epsilon)\tilde x_j(t-T)+\epsilon x_j(t-T)
    \mbox{\quad ($j=1,2$).}
  \end{split}
\end{equation}
We can state two simple corollaries from our general
considerations. First, it is impossible to stabilise the periodic
orbit $x_*$ with extended time-delayed feedback using time-independent
gains $K_1$ and $K_2$ for small $\epsilon$:
\begin{lemma}[Lack of stabilisability for constant control gains]\label{thm:hopf:const}
  Let $p<0$ and let $K_1(t)$ and $K_2(t)$ be arbitrary constants
  (also calling them $K_1$ and $K_2$). Then there exists an
  $\epsilon_{\max}>0$ such that the periodic orbit $x_*$ is unstable
  with the extended time-delayed feedback control \eqref{eq:hopf:etdf}
  for all $\epsilon\in(0,\epsilon_{\max})$.
\end{lemma}
\paragraph{Proof}
Amann \& Hooton \cite{AH13} proved a general topological restriction
on the gains $K(t)$ for extended time-delayed feedback control: let
$K(t)\in\R^n$ be arbitrary (continuous), $\theta\in[0,1]$ be
arbitrary, and let $u$ be of the form
\begin{align*}
  u=\theta K(t)^T[\tilde x(t)-x(t)].
\end{align*}
The scalar $\theta$ provides a homotopy from the uncontrolled system
($\theta=0$) to the controlled system ($\theta=1$). Assume that the
trivial Floquet multiplier $\lambda_1=1$ of $x_*$ is isolated for
$u=0$ (which is the case for example \eqref{eq:hopfu} with $p<0$).
Then the Floquet multiplier $\lambda_1$ depends smoothly on $\theta$
at least for small $\theta$ and will be real: $\lambda_1(\theta)\in\R$
for $0<\theta\ll1$. A necessary condition for extended time-delayed
feedback with gains $K(t)$ to be stabilising for $x_*$ and arbitrary
$\epsilon\in(0,1)$ is that $\lambda_1'(0)\geq0$ if the number of
Floquet multipliers in $\{z\in\C:\Re z>1\}$ is odd for $\theta=0$. If
we denote an adjoint eigenvector for the trivial Floquet multiplier by
$\bar x_*(t)$ (the right eigenvector is $\dot x_*(t)$), this criterion
can be simplified to 
\begin{align*}
  \frac{\int_0^T \bar x_*(t)^Tb(t)K(t)^T\dot x_*(t)\d t}{\int_0^T \bar
    x_*(t)^T\dot x_*(t)\d t}\leq0\mbox{,}
\end{align*}
where $b(t)=\partial_u f(x_*(t),0)$ (this simplifying criterion was
formulated in general in \cite{PN13}). For our particular example, we
have
\begin{align*}
  \dot x_*(t)&=\bar x_*(t)=
  \begin{bmatrix} r\cos t\\ 
    r\sin t
  \end{bmatrix}
\mbox{,} &
  b(t)&=
  \begin{bmatrix}
    1\\ 1
  \end{bmatrix}\mbox{,} &
  T&=2\pi
\end{align*}
and constant gains $K_1$ and $K_2$ such that the
necessary condition of \cite{PN13,AH13} is
\begin{equation}\label{eq:hopf:crit}
  K_1+K_2\leq0\mbox{.}
\end{equation}
On the other hand, if $K_1+K_2\leq0$ the Jacobian of \eqref{eq:hopfu}
with classical linear feedback control 
\begin{align}\label{eq:hopf:classicu}
  u=K_1(x_{1,*}(t)-x_1)+K_2(x_{2,*}(t)-x_2)=K_1(r\sin t-x_1)+K_2(-r\cos t-x_2)
\end{align}
along $x=x_*(t)$ has the trace 
\begin{align*}
  \operatorname{tr}\partial_xf(x(t),K(x_*(t)-x(t)))\vert_{x(t)=x_*(t)}
  =2p+4r^2-K_1-K_2=-2p-K_1-K_2\mbox{.}
\end{align*}
Since $p<0$, this trace is positive if $K_1+K_2\leq0$ for all
$t\in[0,2\pi]$ such that the classical linear feedback control
\eqref{eq:hopf:classicu} cannot be stabilising for the periodic orbit
$x_*=[r\sin t,-r\cos t]^T$. Thus, Lemma~\ref{thm:etdf2classical}
implies that extended time-delayed feedback cannot be stabilising
either, for sufficiently small $\epsilon>0$.

\eop{Lemma~\ref{thm:hopf:const}}

\paragraph{Construction of gains}
\begin{figure}[ht]
  \centering
    \includegraphics[scale=0.5]{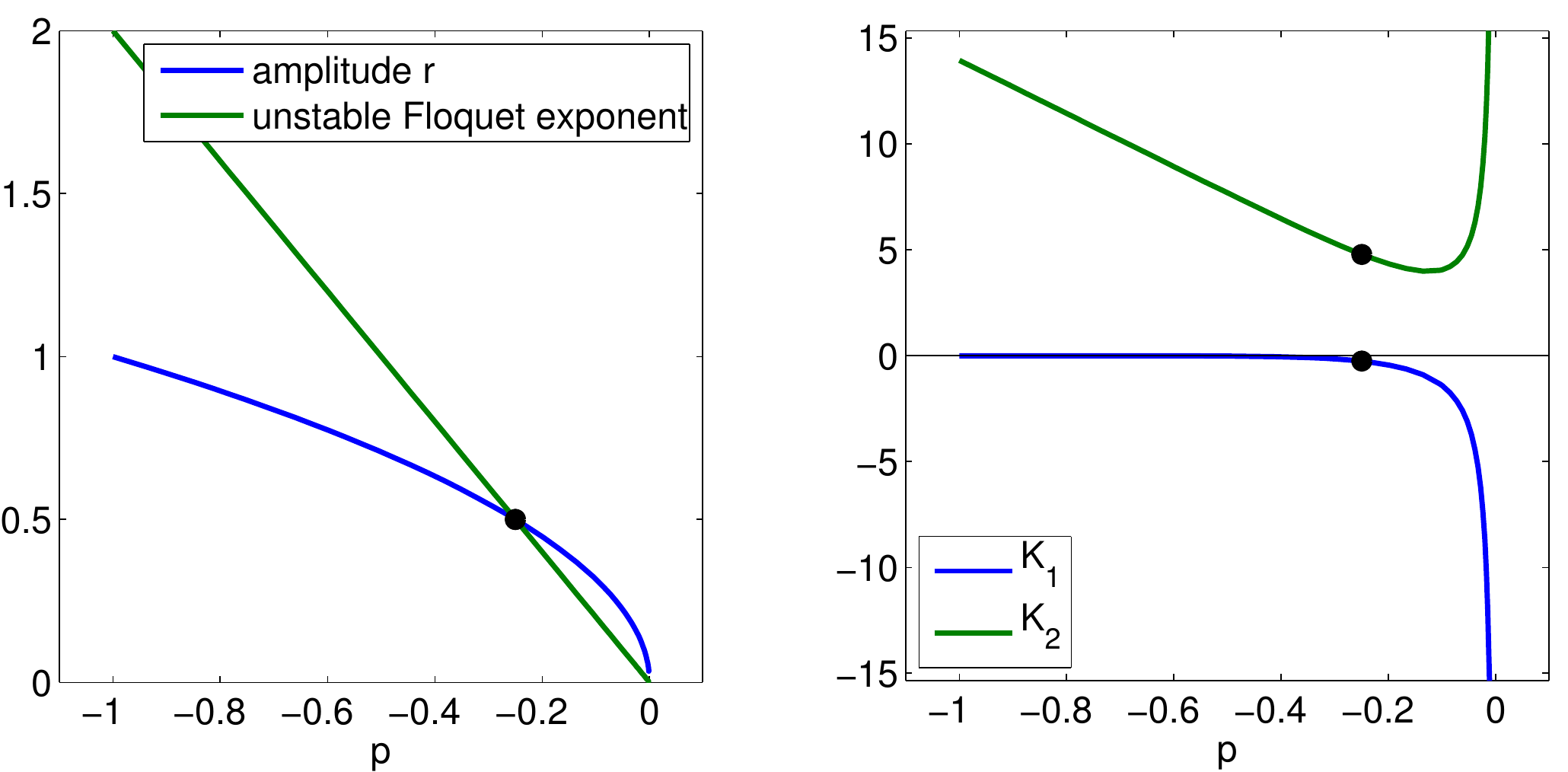}
    \caption{Left: Amplitude $r=\sqrt{-p}$ and unstable Floquet
      exponent (equals $-2p$) along family of periodic orbits. Right:
      gains $K_1$ (note that $K_1<0$ always) and $K_2$ along family of
      periodic orbits in Hopf normal form \eqref{eq:hopfu} with
      feedback law \eqref{eq:hopf:auto:etdf}.}
  \label{fig:gains}
\end{figure}
For the periodic control gains $\Delta_\delta(t)K_0^T$ (or the
autonomous nonlinear gains $\Delta_{\delta,\rho}(x(t))K_0^T$)
the gains $K_0$ are constructed such the matrix $P_0\exp(-b(0)K_0^T)$
has all eigenvalues inside the unit circle (for our illustration we
choose the target location at $\pm\i/2$). Figure~\ref{fig:gains} shows
the amplitude and unstable Floquet exponent of the periodic orbits and
the gains obtained in this manner (called $K_1$ and $K_2$ in
Figure~\ref{fig:gains}). Since the pair $(P_0,b(0))$ is not
controllable at the Hopf point, the gains diverge to infinity for
$p\to0$. In particular, $P_0=\id$ for $p=0$ such that it cannot be
linearly controllable with a single input.
\paragraph{Illustration of asymptotics}
\begin{figure}[ht]
  \centering
  \includegraphics[scale=0.5]{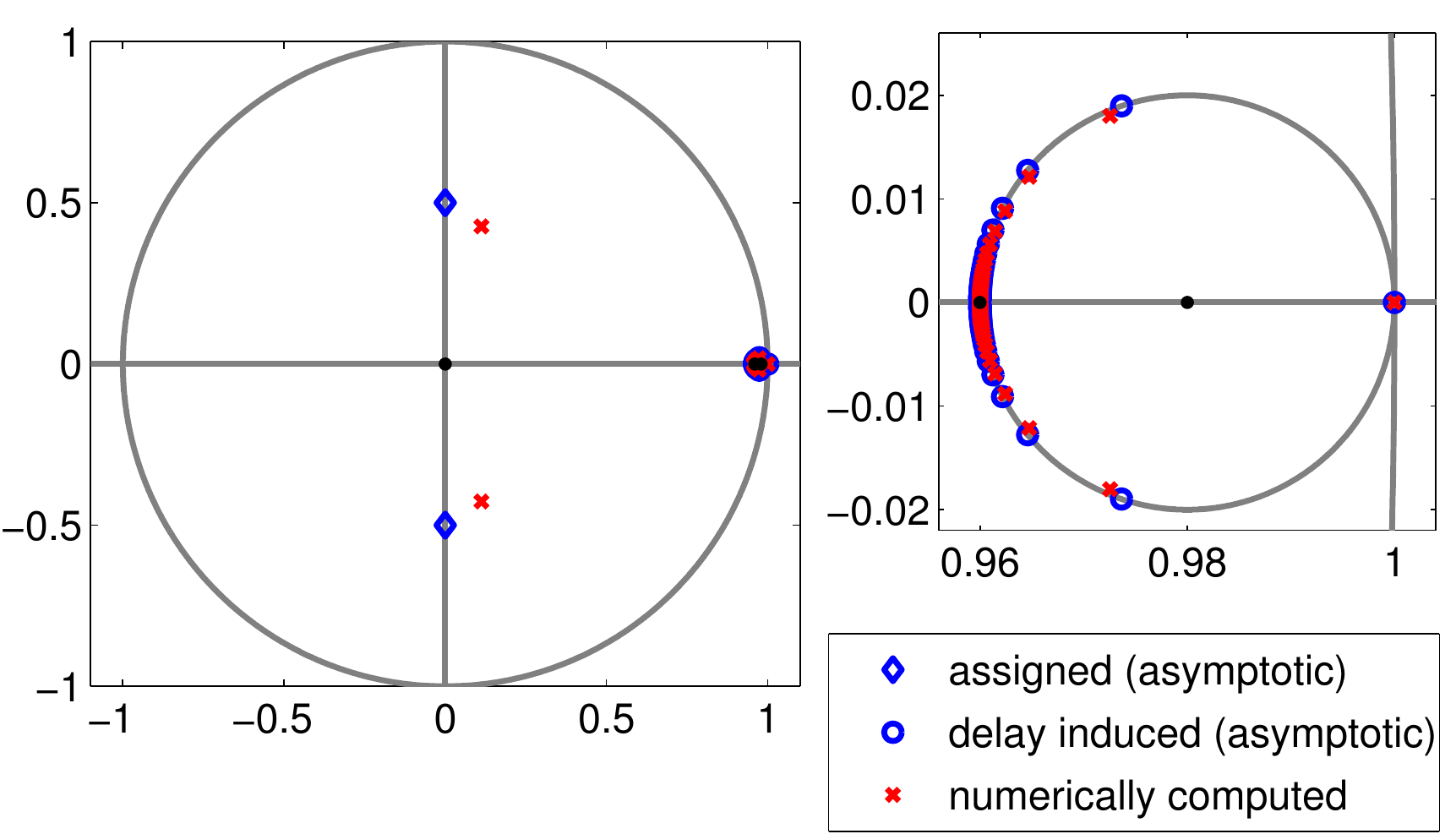}
  \caption{Numerically computed versus asymptotic spectrum for
    $p=-0.25$ ($\epsilon=0.04$, $\delta=T/500$, $\rho=0.3$,
    $K_1=-0.258$, $K_2=4.786$).Left: unit circle in the complex
    plane. Right: zoom into the circle around $1-\epsilon$ or radius
    $\epsilon/2$. Computed with DDE-Biftool
    \cite{ELR02,ELS01,ddebiftoolmanual}, see supplementary material
    and \cite{S16supp} for the code.}
  \label{fig:numeric}
\end{figure}
Figure~\ref{fig:numeric} shows how the true Floquet multipliers
approximate their asymptotic values when using the autonomous
time-delayed feedback control \eqref{eq:etdfnlin:ref} with gains
depending on $x$:
\begin{equation}\label{eq:hopf:auto:etdf}
  \begin{split}
    u(t)&=\Delta_{\delta,\rho}(x(t))
    \left[K_1[\tilde x_1(t)-x_1(t)]+K_2[\tilde x_2(t)-x_2(t)]\right]\\
    \tilde x_j(t)&=(1-\epsilon)\tilde x_j(t-T)+\epsilon x_j(t-T)
    \mbox{\quad ($j=1,2$).}
  \end{split}
\end{equation}
For the construction of $\Delta_{\delta,\rho}$ we used the
construction of $\tilde t$ proposed in
\eqref{eq:auto:u}--\eqref{eq:auto:Jt}:
\begin{align}\label{eq:hopf:trho}
  \tilde t(x)&=\frac{y_0^T}{y_0^Ty_0}[x-x_0]&& \mbox{for $x\in
    B_\rho(x_0)$,}
\end{align}
where $y_0=[r,0]^T$, $x_0=[0,-r]^T$ and $r=\sqrt{-p}$. At the
particular parameter value $p=-0.25$ shown in Figure~\ref{fig:numeric}
the uncontrolled periodic orbit is already strongly unstable: the
unstable Floquet multiplier equals $\exp(-4\pi p)\approx 23.141$. The
gains $K_1$ and $K_2$, designed to assign the Floquet mulitpliers
$\pm\i/2$, are also large after division by $\delta$. Hence, the range
of $\delta$ and $\epsilon$, for which stabilisation is successful is
small. The values for $\delta$ and $\epsilon$ used in the illustration
are chosen such that deviations from the asymptotic limit are visible
but small.

\section{Conclusion and Outlook}
\label{sec:conc}
The paper proves conclusively that there are no restrictions inherent
in extended time-delayed (state) feedback control if one accepts
time-periodic gains, while for constant gains general positive results
are unlikely (as they are absent for classical feedback control of
linear time-periodic systems). The particulars of the gain construction
presented here, following the approach of Brunovksy, are merely for
the purpose of proving their existence analytically. While the result
raises the possibility that more general assignment is feasible (since
there is a lot of freedom in the choice of general time-periodic
$K(t)$) the techniques for proving this may have to be different from
those in the paper. The central argument of the paper rests on the
rank-one nature of the control input, making it easy to locate all
roots $\kappa$ of the transcendental function
$\kappa\mapsto\det[\id-P_0\exp(b(0)K_0^T(\kappa-1))]$ (the matrix $b(0)K_0^T$ has
rank one in our case). The argument also exploits the presence of the
trivial Floquet multiplier, thus, making the result (even if it is
based entirely on linear theory) mostly relevant to the analysis of
nonlinear systems.

\appendix

\section{Auxiliary Lemmas and detailed arguments of proofs}
\label{sec:aux}

\begin{lemma}\label{thm:specass}
  Let $(A,b)$ with $A\in\R^{n\times n}$ and $b\in\R^n$ be
  controllable, let $\det A>0$, and let
  $(\lambda_1,\ldots\lambda_n)\in\C^n$ be the spectrum of a real
  matrix with positive determinant. Then one can find a vector
  $K\in\R^n$ such that
  $\spec(A\exp(bK^T))=(\lambda_1,\ldots,\lambda_n)$.
\end{lemma}
The proof is given in \cite{B69a}. A Matlab implementation of the
explicit construction is \texttt{SpecExpAssign.m} in \cite{S16supp}.

\rev{
\subsection{Floquet multipliers for
  extended time-delayed feedback}
\label{sec:proof:floquet}
\begin{lemma}[Characteristic equation for Floquet
  multipliers]\label{thm:floquet} Let $A(t)\in\R^{n\times n}$, $b(t)$,
  $K(t)\in\R^n$ be $T$-periodic and let $\epsilon$ be positive. Then
  the Floquet multipliers different from $1-\epsilon$ of the linear system
  \begin{align}
    \dot x(t)&=A(t)x(t)+b(t)K(t)^T[\tilde x(t)-x(t)]\label{eq:app:linde}\\
    \tilde x(t)&=(1-\epsilon)\tilde x(t-T)+\epsilon x(t-T)\label{eq:app:linrec}
  \end{align}
  are roots of the function 
  \begin{align*}
    h: \lambda\mapsto \det\left[\lambda\id-P\left(1-
        \frac{\epsilon}{\lambda-(1-\epsilon)}\right)\right]\mbox{.}
  \end{align*}
\end{lemma}
For $\mu\in\C$ the matrix $P(\mu)$ was defined as the solution $x$ at
time $T$ of $\dot x=[A(t)-\mu b(t)K(t)^T]x$ with $x(0)=\id$ (that is,
$P(\mu)$ is the monodromy matrix of $\dot x=[A(t)-\mu b(t)K(t)^T]x$).
\paragraph{Proof} Let $\lambda\in\C$ be an eigenvalue of the time-$T$
map $M$ of \eqref{eq:app:linde}--\eqref{eq:app:linrec}. Let $x_0$,
$\tilde x_0$ (both $[-T,0]\mapsto\R^n$) be the components of an
eigenvector corresponding to $\lambda$, and let $x_1$, $\tilde x_1$
(also both $[-T,0]\mapsto\R^n$) be the corresponding components of
$M[x,\tilde x]$. Then $\tilde x_1(t)=\lambda \tilde x_0(t)$ and, by
definition of the time-$T$ map $M$, $\tilde x_1(t)=(1-\epsilon)\tilde
x_0(t)+\epsilon x_0(t)$. Hence, $\lambda \tilde
x_0(t)=(1-\epsilon)\tilde x_0(t)+\epsilon x_0(t)$, which implies
(since $\lambda\neq1-\epsilon$)
\begin{displaymath}
  \tilde x_0(t)=\frac{\epsilon}{\lambda-(1-\epsilon)} x_0(t)
\end{displaymath}
for $t\in[-T,0]$. Using this relation, we can solve \eqref{eq:app:linde} on the interval $[-T,0]$ as
\begin{displaymath}
  \dot x_0(t)=A(t)x_0(t)+b(t)K(t)^T\left[\frac{\epsilon}{\lambda-(1-\epsilon)}-1\right]x_0(t)
\end{displaymath}
with boundary condition $x_0(0)=\lambda x_0(-T)$. By definition of the
monodromy matrix $P$ this is equivalent to
\begin{displaymath}
  P\left(1-\frac{\epsilon}{\lambda-(1-\epsilon)}\right)x_0(-T)=\lambda x_0(-T)\mbox{,}
\end{displaymath}
which has a non-trivial solution $x_0(-T)$ if and only if the function
$h$ in Lemma~\ref{thm:floquet} is non-zero. \eop{Lemma~\ref{thm:floquet}}
}
\subsection{Proof of Lemma~\ref{thm:etdf:stable}}
\label{sec:proof:etdf:stable}

The characteristic function
$h(\lambda;\epsilon,\delta)=\det[\lambda\id-P(1-\epsilon/(\lambda-(1-\epsilon));\delta,K_0)]$,
defined in \eqref{eq:hdef}, has a root $\lambda=1$ for all small
$\delta$ and all $\epsilon\in(0,1)$: a nullvector of
$\id-P(0;\delta,K_0)$ is $\dot x_*(0)$ (corresponding to a linearised
phase shift).

For $\lambda$ with modulus larger than $1-\epsilon/2$ the term
$\epsilon/(\lambda-(1-\epsilon))$ has modulus less or equal than
$2$. Let us pick $\delta_1>0$ and a $C_1\in(0,1)$ (both small) such that the
polynomial
\begin{align*}
  \lambda\mapsto\det(\lambda\id-P(1-\mu;\delta,K_0))
\end{align*}
has all roots inside the ball $B_{(r+1)/2}(0)\subset B_1(0)$ for all
$\delta\in[0,\delta_1]$ and $\mu$ with $|\mu|\leq C_1$. This is
possible since $\lambda\mapsto \det(\lambda\id-P(1;0,K_0))$ has all
roots inside the ball $B_r(0)$ by assumption of the Lemma and the
limit of $P(1-\mu;\delta,K_0)$ for $\delta\to0$ was uniform for
bounded $\mu$ (recall
$\lim_{\delta\to0}P(1-\mu;\delta,K_0)=P_0\exp((\mu-1) b(0)K_0^T)$).
Thus, $h(\lambda;\epsilon,\delta)$ cannot have roots $\lambda$ on or
outside the unit circle for which
\begin{align*}
  \left|\frac{\epsilon}{\lambda-(1-\epsilon)}\right|\leq C_1
\end{align*}
holds. Hence, for all $\delta\in[0,\delta_1]$ all roots of
$h(\cdot;\epsilon,\delta)$ on or outside of the unit circle must
satisfy
\begin{align}\label{eq:lambda:annulus}
  C_1\leq  \left|\frac{\epsilon}{\lambda-(1-\epsilon)}\right|\leq 2\mbox{.}
\end{align}
We introduce the new variable $\kappa\in\C$ defined via
\begin{align}\label{eq:lambda2kappa}
  \lambda&=1-\epsilon+\epsilon/\kappa\mbox{.}
\end{align}
Restriction \eqref{eq:lambda:annulus} for $\lambda$ is equivalent to the
restriction $ C_1\leq |\kappa|\leq 2$ for $\kappa$.  Hence, for all
$\delta\in[0,\delta_1]$ and $\epsilon\in(0,1)$, every root $\lambda$
of $h(\cdot;\delta,\epsilon)$ on or outside of the unit circle
corresponds to a root $\kappa$ of
\begin{align*}
  g(\kappa;\delta,\epsilon)&=\ h(1-\epsilon+\epsilon/\kappa;\delta,\epsilon)\\
  &=\
  \det\left[\left(1-\epsilon+\frac{\epsilon}{\kappa}\right)\id-
    P(1-\kappa;\delta,K_0)\right]
\end{align*}
with $C_1\leq |\kappa|\leq2$. This one-to-one correspondence of roots
of $h$ and $g$ is given via relation \eqref{eq:lambda2kappa} and
includes multiplicity of the roots. Relation \eqref{eq:lambda2kappa}
also implies that $|\kappa|\leq 1$, because, otherwise, $|\lambda|<1$.
The function $g$ has a limit
\begin{align*}
  g(\kappa;\delta,\epsilon)&\to g(\kappa;0,0)
  &\mbox{for\ }(\delta,\epsilon)\to0
\end{align*}
uniformly for $\kappa$ with $C_1\leq |\kappa|\leq 2$. Hence, the set
of roots $\kappa$ of $g(\cdot;\delta,\epsilon)$ with $C_1\leq
|\kappa|\leq 2$ is a small perturbation of the set of roots of $g(\cdot;0,0)$:
\begin{align*}
  g(\kappa;0,0)&=\det\left[\id-
    P_0\exp\left(b(0)K_0^T
      \left[\kappa-1\right]\right)\right]\\
  &=\det\left[\id-
    P_0-P_0b(0)K_0^T\sigma(\kappa-1)\right]\mbox{, where }\\
  \sigma(\kappa-1)&=
  \begin{cases}
    \cfrac{\exp\left(K_0^Tb(0)(\kappa-1)\right)-1}{K_0^Tb(0)} 
    &\mbox{if $K_0^Tb(0)\neq0$}\\
    \kappa-1 &\mbox{if $K_0^Tb(0)=0$.}
  \end{cases}
\end{align*}
The generalised eigenvalue problem for the matrix pair
$(\id-P_0,P_0b(0)K_0^T)$ with characteristic polynomial
$\sigma\mapsto\det(\id-P_0-P_0b(0)K_0^T\sigma)$ is regular because
$\det(\id-P_0)=0$ but $\det(\id-P_0-P_0b(0)K_0^T\sigma(-1))$ is
regular (since $P_0\exp(-b(0)K_0^T)$ has all eigenvalues inside the
unit circle, $\id-P_0\exp(-b(0)K_0^T)$ is regular). As
$P_0b(0)K_0^T$ has rank $1$ the characteristic polynomial
corresponding to $\sigma\mapsto \det(\id-P_0-P_0b(0)K_0^T\sigma)$ has
degree $1$. Moreover, its only root equals $0$ (which must be simple
due to the regularity of $\det(\id-P_0-P_0b(0)K_0^T\sigma)$). Hence,
we know that $g(\kappa;0,0)=0$ if and only if $\sigma(\kappa-1)=0$.

\paragraph{Case $K_0^Tb(0)=0$} If $K_0^Tb(0)=0$, this implies that the
only root $\kappa$ of $g(\cdot;0,0)$ with $C_1\leq|\kappa|\leq 2$
equals unity. Hence, also for sufficiently small $\epsilon$ and
$\delta$, the only root $\kappa$ with $C_1\leq|\kappa|\leq 2$ of
$g(\cdot;\delta;\epsilon)$ equals unity (since $g(1;\delta;\epsilon)=0$).

\paragraph{Case $K_0^Tb(0)\neq0$} If $K_0^Tb(0)\neq0$, we have that
$g(\kappa;0,0)=0$ if and only if $\exp(K_0^Tb(0)(\kappa-1))=1$ such that
the roots are
\begin{align*}
  \kappa_{\ell,0}&=1+\frac{2\pi\i\ell}{K_0^Tb(0)} &&\mbox{($\ell\in\Z$)}
\end{align*}
The first part of the subscript, $\ell$, numbers the roots, the second part of the subscript, $0$, indicates that
$\delta=\epsilon=0$.  Thus, the roots $\kappa_{\ell,0}$ of
$g(\cdot;0,0)$ have a modulus
\begin{align*}
  |\kappa_{\ell,0}|=\sqrt{1+\frac{4\pi^2\ell^2}{(K_0^Tb(0))^2}}
\end{align*}
such that only the roots $\kappa_{\ell,0}$ with index
\begin{align*}
  -\ell_{\max}&\leq\ell\leq\ell_{\max}\mbox{,} &
  \mbox{where\ }\ell_{\max}&=\frac{\sqrt{3}}{2\pi}|K_0^Tb(0)|\mbox{,}
\end{align*}
are in the admissible range with $C_1\leq|\kappa_{\ell,0}|\leq
2$. (Hence, both cases, $K_0^Tb(0)=0$ and $K_0^Tb(0)\neq0$ can be
treated equally.) The admissible roots $\kappa_{\ell,0}$ of
$g(\cdot;0,0)$ are all simple. Hence, for sufficiently small $\delta$
and $\epsilon$, $g(\cdot;\delta,\epsilon)$ will have roots
$\kappa_\ell$ for $|\ell|\leq\ell_{\max}$ that are small perturbations
of $\kappa_{\ell,0}$, and these roots $\kappa_\ell$ are the only roots
of $g(\cdot;\delta,\epsilon)$ with modulus in $[C_1,2]$. Since we know
that $g(1;\delta,\epsilon)=0$, we know that $\kappa_0=1$ (hence, for
$\ell=0$ the perturbation is zero). Furthermore, for non-zero $\ell$
with $|\ell|\leq\ell_{\max}$, the modulus of $\kappa_{\ell,0}$ is
greater than $1$. Hence, the perturbed roots $\kappa_\ell$ also have
modulus greater then $1$ for sufficiently small $\epsilon$ and
$\delta$, and non-zero $|\ell|\leq\ell_{\max}$.

Consequently, by relation~\eqref{eq:lambda2kappa}, the only roots of
$h(\cdot;\delta,\epsilon)$ that could be on or outside the unit circle
are
\begin{align*}
  \lambda_\ell&=1-\epsilon+\frac{\epsilon}{\kappa_\ell} &
  \mbox{where\ } |\ell|\leq\ell_{\max}\mbox{.}
\end{align*}
However, these roots $\lambda_\ell$ are simple and satisfy
$\lambda_0=1$ and $|\lambda_\ell|<1$ for non-zero $\ell$, since
$|\kappa_\ell|>1$ for non-zero $\ell$.
\eop{Lemma~\ref{thm:etdf:stable}}

\subsection{Details of construction for autonomous feedback gains --- Regularisation of the short impulse $\Delta_\delta(t)$}
\label{sec:ktmod}
To avoid discontinous dependence of the right-hand side on the
solution, we first regularise the discontinuity of the time-dependent
gain $K(t)$. Define for $\delta\in(0,\sqrt{T+1/16}-1/4)$ (such that
$2\delta^2+\delta<T$) the regularised version of $\Delta_\delta$:
\begin{align}\label{eq:app:delta}
  \Delta_\delta(t)&=
  \begin{cases}
        1/\delta & \mbox{if $t\vert_{\mod[0,T)}\in[0,\delta]$,}\\
        0 & \mbox{if $t\vert_{\mod[0,T)}\in[\delta+\delta^2,T-\delta^2]$,}\\
        \cfrac{1}{\delta}\ m\left(\cfrac{\delta+\delta^2-t\vert_{\mod[0,T)}}{\delta^2}\right)
        &\mbox{if $t\vert_{\mod[0,T)}\in(\delta,\delta+\delta^2)$}\\[2.5ex]
        \cfrac{1}{\delta}\ m\left(\cfrac{t\vert_{\mod[0,T)}-T+\delta^2}{\delta^2}\right)
        &\mbox{if $t\vert_{\mod[0,T)}\in(T-\delta^2,T)$,}        
  \end{cases}
\end{align}
where $m:\R\mapsto[0,1]$ is an arbitrary smooth monotone increasing
function with $m(s)=0$ for $s\leq0$ and $m(s)=1$ for $s\geq1$. When
using $\Delta_\delta$ as defined in \eqref{eq:app:delta} instead of
\eqref{eq:intro:delta} to define the linear (now approximately)
short-impulse feedback law
\begin{align}\label{eq:utildedef}
  u_\delta(t;y)&=\Delta_\delta(t)K_0^Ty
\end{align}
the nonlinear time-$T$ map is still linearisable. Denoting the
monodromy matrix of the linear system (recall that
$A(t)=\partial_xf(x_*(t),0)$, $b(t)=\partial_uf(x_*(t),0)$ and using
definition \eqref{eq:app:delta} for $\Delta_\delta$)
\begin{align*}
  \dot y(t)&=\left[A(t)-\mu b(t)
    \Delta_\delta(t)K_0^T\right]y(t)\mbox{,}& y(0)&=y_0
\end{align*}
again by $P(\mu;\delta,K_0)$ then the monodromy matrix of the smoothed
system still satisfies (identical to \eqref{eq:Plimit})
\begin{align}
  \label{eq:ptlimit}
  P(\mu;\delta,K_0)= P_0\exp(-b(0)K_0^T\mu)+O(\delta)\mbox{,}
\end{align}
where the error term $O(\delta)$ is uniform for bounded $\mu\in\C$ and
$K_0\in\R^n$. Hence, we can replace the discontinuous definition
\eqref{eq:intro:delta} for $\Delta_\delta(t)$ by \eqref{eq:app:delta}
in \eqref{eq:etdfode}, and Lemma~\ref{thm:etdf:stable} still applies
to the modified (regularised) system.

\subsection{State-dependent gains $K(x(t))$}
\label{sec:sdK}
Consider a sufficiently small radius $\rho>0$ such that the equation
$\dot x_*(0)^T[x-x_*(t)]=0$ has a unique solution $t\in\R$ close to
$0$ for all $x\in B_{2\rho}(x_*(0))\subset\R^n$, thus defining
implicitly a smooth function $t_\rho$:
\begin{align}\label{eq:tx:implicit}
  t_\rho:&\ \R^n\mapsto \R\mbox{,}& t_\rho(x)&=
  \begin{cases}
    \mbox{root $t$ of $\dot x_*(0)^T[x-x_*(t)]=0$} &
    \mbox{if $x\in B_{2\rho}(x_*(0))$,}\\
    \mbox{arbitrary such that $ t_\rho$ is smooth} &
    \mbox{otherwise.}
  \end{cases}
\end{align}
The function $\tilde t$, defined in \eqref{eq:auto:Jt}, is
approximately equal to $t_\rho$ along the periodic orbit $x_*$ and
near $x_*(0)$: $t_\rho(x_*(t))-\tilde t(x_*(t))=t-\tilde
t(x_*(t))=O(t^2)$. We consider also a regularised indicator function
for the neighbourhood of $x_*(0)$
\begin{align*}
  J_\rho:&\ \R^n\mapsto\R\mbox{,} &
  J_\rho&=\begin{cases}
    1 &
    \mbox{if $x\in B_\rho(x_*(0))$,}\\
    0 &     \mbox{if $x\notin B_{2\rho}(x_*(0))$,}\\
    \mbox{arbitrary such that $J_\rho$ is smooth} &
    \mbox{if $x\in B_{2\rho}(x_*(0))\setminus B_\rho(x_*(0))$.}
  \end{cases}
\end{align*}
and combine $ t_\rho$ and $J_\rho$ with $\tilde \Delta_\delta$ as
defined in \eqref{eq:app:delta} to the smooth globally defined
function
\begin{align*}
  \tilde \Delta_{\delta,\rho}:&\ \R^n\mapsto \R\mbox{,} &
  \tilde \Delta_{\delta,\rho}(x)&=J_\rho(x) \tilde \Delta_\delta( t_\rho(x))
\end{align*}
When applying $\tilde \Delta_{\delta,\rho}$ to $x_*(t)$ the result is
identical to the timed impulse $\tilde \Delta_\delta$ for small
$\delta$: if $\|x_*(t)-x_*(0)\|<\rho$ for all
$t\in[-\delta^2,\delta+\delta^2]$ then
\begin{align*}
  \tilde \Delta_{\delta,\rho}(x_*(t))=\tilde \Delta_\delta(t)
\end{align*}
for all $t\in\R$. Consequently, the system with extended time-delayed
feedback and state-dependent gains
\begin{align}
  \label{eq:etdfnlin:ode}
  \dot x(t)&=f(x(t),\tilde \Delta_{\delta,\rho}(x(t))K_0^T[\tilde x(t)-x(t)])\\
  \tilde x(t)&=(1-\epsilon)\tilde x(t-T)+\epsilon x(t-T)\mbox{,}
  \label{eq:etdfnlin:ref}
\end{align}
which is now autonomous with a smooth right-hand side, has for
sufficiently small $\rho$ and $\delta+\delta^2<\rho$ exactly the same
linearisation along the periodic orbit $x(t)=\tilde x(t)=x_*(t)$ as
system~\eqref{eq:etdfode},\,\eqref{eq:etdfref}. The derivative of
$\tilde \Delta_{\delta,\rho}$ with respect to its argument is multiplied by
$0$ if $\tilde x(t)=x(t)=x_*(t)$ for all times $t$ in the term $
\tilde \Delta_{\delta,\rho}(x(t))K_0^T[\tilde x(t)-x(t)]$ in
\eqref{eq:etdfnlin:ode}.

The time reconstruction function $t_\rho$ as defined in
\eqref{eq:tx:implicit} satisfies $t_\rho(x_*(t))=t$ as long as
$x_*(t)\in B_{2\rho}(x_*(0))$. The definition of $t_\rho$ as proposed
in \eqref{eq:auto:Jt} in the main text is a $O(t^2)$ perturbation of
\eqref{eq:tx:implicit} for $t$ of order $\delta$. Thus, continuity of
the Flqouet multipliers implies that the perturbed version of $t_\rho$
preserves the stability of the periodic orbit $x_*$.

\subsubsection*{Declaration}
\label{declaration}
\begin{description}
\item[Competing interests] I have no competing interests.
\item[Access to data and code] on 
  \url{https://dx.doi.org/10.6084/m9.figshare.2993812}.
\item[Funding] Jan Sieber's research has received funding from the
  European Union's Horizon 2020 research and innovation programme
  under Grant Agreement number 643073.
\end{description}


\begin{thebibliography}{10}

\bibitem{P92}
Pyragas K.
\newblock {Continuous control of chaos by self-controlling feedback}.
\newblock Phys Lett A. 1992;170:421--428.

\bibitem{PT93a}
Pyragas K, Tama\v{s}evi\v{c}ius A.
\newblock Experimental control of chaos by delayed self-controlling feedback.
\newblock Physics Letters A. 1993;180(1-2):99--102.

\bibitem{SSG94}
Socolar JES, Sukow DW, Gauthier DJ.
\newblock {Stabilizing unstable periodic orbits in fast dynamical systems}.
\newblock Phys Rev E. 1994;50:3245.

\bibitem{KBPOMBRE01}
Kim M, Bertram M, Pollmann M, Oertzen Av, Mikhailov AS, Rotermund HH, et~al.
\newblock {Controlling chemical turbulence by global delayed feedback{\rm:}
  pattern formation in catalytic {CO} oxidation on {P}t(110)}.
\newblock Science. 2001;292(5520):1357--1360.

\bibitem{UBBWH04}
Ushakov O, Bauer S, Brox O, W\"{u}nsche HJ, Henneberger F.
\newblock {Self-Organization in Semiconductor Lasers with Ultra-Short Optical
  Feedback}.
\newblock Phys Rev Lett. 2004;92:043902.

\bibitem{WSH08}
W{\"u}nsche HJ, Schikora S, Henneberger F.
\newblock Noninvasive control of semiconductor lasers by delayed optical
  feedback.
\newblock In: Schuster HG, Sch{\"o}ll E, editors. Handbook of Chaos Control.
  2nd ed. Wiley-VCH; 2008. p. 455--474.

\bibitem{SHWSH06}
Schikora S, H\"{o}vel P, W\"{u}nsche HJ, Sch\"{o}ll E, Henneberger F.
\newblock {All-Optical Noninvasive Control of Unstable Steady States in a
  Semiconductor Laser}.
\newblock Phys Rev Lett. 2006;97(21):213902.

\bibitem{PHT05}
Popovych OV, Hauptmann C, Tass PA.
\newblock Effective desynchronization by nonlinear delayed feedback.
\newblock Physical Review letters. 2005;94(16):164102.

\bibitem{PHT06}
Popovych OV, Hauptmann C, Tass PA.
\newblock Control of neuronal synchrony by nonlinear delayed feedback.
\newblock Biological Cybernetics. 2006;95(1):69--85.

\bibitem{SHHD09}
Sch{\"o}ll E, Hiller G, H{\"o}vel P, Dahlem MA.
\newblock Time-delayed feedback in neurosystems.
\newblock Philosophical Transactions of the Royal Society of London A:
  Mathematical, Physical and Engineering Sciences. 2009;367(1891):1079--1096.

\bibitem{LWP01}
L\"uthje O, Wolff S, Pfister G.
\newblock Control of Chaotic Taylor-Couette Flow with Time-Delayed Feedback.
\newblock Phys Rev Lett. 2001;86:1745--1748.

\bibitem{YH06}
Yamasue K, Hikihara T.
\newblock Control of microcantilevers in dynamic force microscopy using time
  delayed feedback.
\newblock Review of scientific instruments. 2006;77(5):053703.

\bibitem{SGNWK08}
Sieber J, Gonzalez-Buelga A, Neild SA, Wagg DJ, Krauskopf B.
\newblock {Experimental continuation of periodic orbits through a fold}.
\newblock Phys Rev Lett. 2008;100:244101.

\bibitem{BMB12}
Barton DAW, Mann BP, Burrow SG.
\newblock Control-based continuation for investigating nonlinear experiments.
\newblock Journal of Vibration and Control. 2012;18(4):509--520.

\bibitem{BSESTS2014}
Bureau E, Schilder F, Elmeg{\aa}rd M, Santos IF, Thomsen JJ, Starke J.
\newblock Experimental bifurcation analysis of an impact oscillator ---
  Determining stability.
\newblock Journal of Sound and Vibration. 2014;333(21):5464--5474.

\bibitem{MSA04}
Montagnier P, Spiteri RJ, Angeles J.
\newblock The control of linear time-periodic systems using
  {F}loquet-{L}yapunov theory.
\newblock International Journal of Control. 2004;77(5):472--490.

\bibitem{B69a}
Brunovsky P.
\newblock Controllability and linear closed-loop controls in linear periodic
  systems.
\newblock Journal of Differential equations. 1969;6(2):296--313.

\bibitem{HL93}
Hale JK, {Verduyn Lunel} SM.
\newblock {Introduction to functional-differential equations}. vol.~99 of
  Applied Mathematical Sciences.
\newblock New York: Springer-Verlag; 1993.

\bibitem{MEVR02}
Michiels W, Engelborghs K, Vansevenant P, Roose D.
\newblock Continuous pole placement for delay equations.
\newblock Automatica. 2002;38(5):747--761.

\bibitem{MN14}
Michiels W, Niculescu SI.
\newblock Stability, Control, and Computation for Time-Delay Systems: An
  Eigenvalue-Based Approach. vol.~27.
\newblock SIAM; 2014.

\bibitem{LHFGFS11}
Lehnert J, H{\"o}vel P, Flunkert V, Guzenko PY, Fradkov AL, Sch{\"o}ll E.
\newblock Adaptive tuning of feedback gain in time-delayed feedback control.
\newblock Chaos: An Interdisciplinary Journal of Nonlinear Science.
  2011;21(4):043111.

\bibitem{PP13}
Pyragas V, Pyragas K.
\newblock Adaptive search for the optimal feedback gain of time-delayed
  feedback controlled systems in the presence of noise.
\newblock The European Physical Journal B. 2013;86(7):1--8.

\bibitem{PP14}
Pyragas V, Pyragas K.
\newblock Continuous pole placement method for time-delayed feedback controlled
  systems.
\newblock The European Physical Journal B. 2014;87(11):1--10.

\bibitem{PP11}
Pyragas V, Pyragas K.
\newblock Adaptive modification of the delayed feedback control algorithm with
  a continuously varying time delay.
\newblock Physics Letters A. 2011;375(44):3866--3871.

\bibitem{NP12}
Novi{\v{c}}enko V, Pyragas K.
\newblock Phase-reduction-theory-based treatment of extended delayed feedback
  control algorithm in the presence of a small time delay mismatch.
\newblock Physical Review E. 2012;86(2):026204.

\bibitem{JBORB97}
Just W, Bernard T, Ostheimer M, Reibold E, Benner H.
\newblock {Mechanism of Time-Delayed Feedback Control}.
\newblock Phys Rev Lett. 1997;78(2):203--206.

\bibitem{SS07}
Sch\"{o}ll E, Schuster HG, editors.
\newblock {Handbook of Chaos Control}.
\newblock 2nd ed. New York: Wiley; 2007.

\bibitem{N04}
Nakajima H.
\newblock {Some sufficient conditions for stabilizing periodic orbits without
  the odd-number property by delayed feedback control in continuous-time
  systems}.
\newblock Phys Lett A. 2004;327:44--54.

\bibitem{NU98}
Nakajima H, Ueda Y.
\newblock {Limitation of generalized delayed feedback control}.
\newblock Physica D. 1998;111:143--150.

\bibitem{P01}
Pyragas K.
\newblock {Control of chaos via an unstable delayed feedback controller}.
\newblock Phys Rev Lett. 2001;86(11):2265--2268.

\bibitem{FFGHS07}
Fiedler B, Flunkert V, Georgi M, H\"{o}vel P, Sch\"{o}ll E.
\newblock {Refuting the Odd-Number Limitation of Time-Delayed Feedback
  Control}.
\newblock Phys Rev Lett. 2007;98(11):114101.

\bibitem{FS11}
Flunkert V, Sch\"oll E.
\newblock Towards easier realization of time-delayed feedback control of
  odd-number orbits.
\newblock Phys Rev E. 2011 Jul;84:016214.

\bibitem{PN13}
Pyragas K, Novi{\v{c}}enko V.
\newblock Time-delayed feedback control design beyond the odd-number
  limitation.
\newblock Physical Review E. 2013;88(1):012903.

\bibitem{BPS11}
Postlethwaite GBCM, Silber M.
\newblock Time-delayed feedback control of unstable periodic orbits near a
  subcritical Hopf bifurcation.
\newblock Physica D: Nonlinear Phenomena. 2011;240(9):859--871.

\bibitem{PBS13}
Postlethwaite CM, Brown G, Silber M.
\newblock Feedback control of unstable periodic orbits in equivariant Hopf
  bifurcation problems.
\newblock Philosophical Transactions of the Royal Society of London A:
  Mathematical, Physical and Engineering Sciences. 2013;371(1999):20120467.

\bibitem{SHFD10}
Sch{\"o}ll E, H{\"o}vel P, Flunkert V, Dahlem MA.
\newblock Time-delayed feedback control: from simple models to lasers and
  neural systems.
\newblock In: Complex time-delay systems. Springer; 2010. p. 85--150.

\bibitem{HA12}
Hooton EW, Amann A.
\newblock Analytical limitation for time-delayed feedback control in autonomous
  systems.
\newblock Physical Review letters. 2012;109(15):154101.

\bibitem{AH13}
Amann A, Hooton EW.
\newblock An odd-number limitation of extended time-delayed feedback control in
  autonomous systems.
\newblock Philosophical Transactions of the Royal Society of London A:
  Mathematical, Physical and Engineering Sciences. 2013;371(1999):20120463.

\bibitem{S16supp}
Sieber J.
\newblock Generic stabilisability for time-delayed feedback control
  (Supplementary material).
\newblock Figshare. 2016;\url{https://dx.doi.org/10.6084/m9.figshare.2993812}.

\bibitem{ELR02}
Engelborghs K, Luzyanina T, Roose D.
\newblock Numerical bifurcation analysis of delay differential equations using
  {DDE-BIFTOOL}.
\newblock ACM Transactions on Mathematical Software. 2002;28(1):1--21.

\bibitem{ELS01}
Engelborghs K, Luzyanina T, Samaey G.
\newblock {{DDE-BIFTOOL} v.2.00{\rm:} a {Matlab} package for bifurcation
  analysis of delay differential equations}.
\newblock Katholieke Universiteit Leuven; 2001. 330.

\bibitem{ddebiftoolmanual}
Sieber J, Engelborghs K, Luzyanina T, Samaey G, Roose D. DDE-BIFTOOL Manual ---
  Bifurcation analysis of delay differential equations;.

\end{thebibliography}

\end{document}